\documentclass[a4paper,12pt]{article}
\usepackage{theorem,amsmath,latexsym,amssymb,newcom,esint}
\usepackage[mathscr]{eucal}
\usepackage[dvips]{graphicx,color}
\usepackage{wrapfig}
\begin{document}
\title{\textsf{A boundary partial regularity \\
               and a regularity criterion for \\
               New Harmonic Heat flows}}
\author {\textsf{Kazuhiro HORIHATA}}
\maketitle
\begin{abstract}
In my previous paper of K.Horihata~\cite{horihata}, 
 we have proposed a Ginzburg-Landau type heat flow 
  with a time-dependent parameter
   and then passing to the limit of the parameter appeared in it, 
    we have constructed 
     a harmonic heat flow into spheres.
By using this scheme, we establish a few energy inequalities:

(i) monotonical inequalities and

(ii) a reverse Poincare inequality at any boundary point.

Furthermore we assert that these inequalities (i) and (ii) derive
 the smaller estimates on the set of noncontinuous points
  for our flow contrast to the former results.
We refer to them by Y.Chen~\cite{chen-91} and Y.Chen, J.Li, F.H.Lin~\cite{chen-li-lin}.
 Next we introduce two conditions for the whole domain's
  regularity for it; 
The one is a boundary energy smallness and the another 
 is an one-sided condition proposed by S.Hildebrandt and K.-O.Widman~\cite{hildebrandt-widman}.

\end{abstract}
%
%
%
\setcounter{chapternumber}{1}\setcounter{equation}{0}
\renewcommand{\theequation}%
           {\thechapternumber.\arabic{equation}}
\section{\enspace Introduction.}
\label{SEC:Intro}
A previous paper: K.Horihata~\cite{horihata}, constructed 
 a harmonic heat flow between $d$-dimensional unit ball and
  $D$-dimensional unit sphere with
\renewcommand{\theenumi}{\roman{enumi}}
\begin{enumerate}
\item
a global energy inequality,
\item
a monotonical inequality
\item
a reverse Poincar\'e inequality
\end{enumerate}
and showed the flow is smooth except a small singular set.
This paper discusses a partial regularity for such a flow 
 from $C^2$-domain in $\mathbb{R}^d$ to
  the $D$-dimensional unit sphere
   and then state two condition for 
    the whole domain's smoothness for it under a certain hypothesis: 
%
For any point $x_0$ $\in$ $\partial \Omega$,
 there exist a neighbourhood $U(x_0)$ and
  $C^2$-diffeomorphism $\Phi_{x_0}$ 
    between $U(x_0)$ and $\mathbb{R}^d$, writing
\allowdisplaybreaks\begin{align}
& 
y \; = \; ( y^\prime, y_d ) \; = \; 
 \Phi_{x_0} (x) \; = \;
  (x^\prime, x_d \, - \, \phi_{x_0} (x^\prime))
\notag
\\
&
\quad ( y \, = \, ( y^\prime, y_d) \, \in \, \mathbb{R}^d,
 \; x \, = \, (x^\prime,x_d)
  \, \in \, U(x_0))
\notag
\end{align}
whereupon writing $x_0$ $=$ $(x_0^\prime, x_{d,0})$, 
 a suitable $C^2$-function $\phi_{x_0}$ 
  on $U(x_0)$ $\cap$ $\mathbb{R}_0^d$ satisfies
\allowdisplaybreaks\begin{align}
&
\phi_{x_0} (x_0^\prime) \; = \; 0,
 \; \nabla \phi_{x_0} (x_0^\prime) \; = \; 0,
\notag
\\
&
\sum_{i,j=1}^{d-1}
 \nabla_{i} \nabla_j \phi_{x_0} ( x_0^\prime ) 
  ( x^\prime - x_0^\prime)_i ( x^\prime - x_0^\prime)_j 
   \; \ge \; 
    \theta_0 \vert x^\prime - x_0^\prime \vert^2
\label{Cond:B}
\tag{B}
\end{align}
for any point $x$ $=$ $(x^\prime, x_d)$ 
 $\in$ $U(x_0)$ $\cap$ $\partial \Omega$,
  where a positive constant $\theta_0$ is independent of $x_0$.
We call the condition our the boundary $(\textit{B})$.
\par
We rigorously state our problem and results:
 Let $\Omega$ and $\mathbb{S}^D$ be respectively 
  the $C^2$-domain in $\mathbb{R}^d$ (\textit{not always bounded}) and
   the unit sphere in $\mathbb{R}^{D+1}$
    and a parabolic cylinder $Q(T)$ be $(0,T)$ $\times$ $\Omega$,
     where $d$ and $D$ are positive integers greater than or equal to $2$.
Consider the Sobolev class
 $H_B^{k,p} (\Omega ;\mathbb{S}^D)$ $\, := \,$
  $\{u \in L (\Omega \cap B_n (0); \mathbb{R}^{D+1})$,
   $\nabla u \in H^{k-1,p}$ $(\Omega ; \mathbb{R}^{D+1})$
    $; |u| = 1 \quad \mathrm{a.e.} \, \, x \in \Omega \}$.
     for any positive integer $n$,
      where $k$ is any nonnegative and $p$ is positive integers
       greater than equal to $1$.
Giving a mapping $u_0 \, \in \, $
 $H_B^{1,2} (\Omega ; \mathbb{S}^D)$,
  we consider the heat flow:
\begin{equation}
\left\{
\begin{array}{rll}
\dfrac {\partial u}{\partial t} & \, = \,
 \triangle u \, + \, | \nabla u |^2 u
  \quad & \mathrm{in} \;\; Q(T),
\\[5pt]
u(0,x) & \, = \, u_0(x)
 \qquad & \mathrm{at} \quad \{0\} \times \Omega,
\\[3pt]
u(t,x) & \, = \, u_0(x)
 \qquad & \mathrm{on} \quad [0,T) \times \partial \Omega.
\end{array}
\right.
\label{EQ:HHF}
\end{equation}
Define three function classes:
\allowdisplaybreaks\begin{align}
L^\infty & \bigl(0,T;H_B^{1,2}(\Omega; \mathbb{S}^D)\bigr) 
\notag
\\
&
\, := \, 
 \{ u | u \; \text{is measurable from} \;
  [0,T] \;\mathrm{to}\; H^{1,2} (\Omega ;\mathbb{S}^D) 
\notag
\\
H^{1,2} & \bigl(0,T; L_B^2(\Omega;\mathbb{R}^{D+1}) \bigr) 
\notag
\\
&
\; := \;
 \Bigl\{ u \, \in \, L_B^2 (\Omega;\mathbb{R}^{D+1}) 
  \, ; \, \partial u/\partial t
   \, \in \, L^2 (Q(T) ; \mathbb{R}^{D+1}) \Bigr\}
\notag
\\
&
V_B (Q(T);\mathbb{S}^D) \; := \; 
 L^\infty \bigl(0,T;H_B^{1,2}(\Omega;\mathbb{S}^D)\bigr) \cap 
  H^{1,2} \bigl(0,T;L_B^2(\Omega;\mathbb{R}^{D+1})\bigr).
\notag
\end{align}
The weak formulation of \eqref{EQ:HHF} is as follows:
 For any given mapping $u_0 \, \in \,$ 
  $H_B^{1,2} (\Omega;\mathbb{S}^D)$,
   we call a mapping $u \, \in \, V_B (Q(T);\mathbb{S}^D)$
    a weakly harmonic heat flow {(\textit{WHHF})} provided
for any $\phi \in C_0^\infty ({Q(T)};\mathbb{R}^{D+1})$
\allowdisplaybreaks\begin{alignat}{2}
&
\lint_{Q(T)} \Bigl( \la \frac {\partial u}{\partial t}, \phi \ra
 \, + \, \la \nabla u, \nabla \phi \ra
  \, - \, \Bigr. && \Bigl. 
   \la u, \phi \ra |\nabla u|^2
    \Bigr) \, dz \, = \,0,
\label{EQ:1}
\\
& 
u(t) \, - \, u_0 \in
 \tc_B (\Omega ;\mathbb{R}^{D+1}) 
  &&\quad \text{for almost every} \;
   t \, \in \, (0,T),
\label{EQ:2}
\\
& \underset {t \searrow 0}{\lim} \,
 u(t) \, = \, u_0
  \quad && \quad \mathrm{in} \quad L_B^2 (\Omega ;\mathbb{R}^{D+1}).
\label{EQ:3}
\end{alignat}
\par
K.Horihata~\cite{horihata} has proposed a new approximate evolutional scheme
 said to be the Ginzburg-Landau heat flow (\textit{GLHF}).
We review Ginzburg-Landau heat flow:
 Introduce smooth functions $\chi (t)$ and $\kappa (t)$ by
\allowdisplaybreaks\begin{align}
&
\chi (t) \; = \;
\begin{cases}
t & \quad (t \, < \, 2)
\\
3 & \quad (t \, \ge \, 4),
\end{cases}
\notag
\\
&
\kappa (t) \; = \; \arctan (t)/ \pi.
\end{align}
Let a mapping $u_0$ in 
 $H_B^{1,2}$ $( \Omega \, ; \, \mathbb{S}^D )$
  $\cap$ 
   $(H_B^{2-1/q,q} \cap H_B^{[(d+1)/2]+1,p})$ 
    $( \Omega\setminus (\complement\Omega)_{\delta_0} ; \mathbb{S}^D)$
     with positive numbers sufficiently small $\delta_0$,
      and positive numbers $q$ greater than or equal to $d+2$ and 
       $p$ greater than $2$.
Then GLHF is designated by solution of the systems:
\begin{eqnarray}
\dfrac{\partial u_\lambda}{\partial t} \, - \,
 \triangle u_\lambda
  \, + \, \lambda^{1-\kappa}
   \dot \chi \bigl( \bigl( | u_\lambda |^2 \, - \, 1 \bigr)^2 \bigr)
    ( | u_\lambda |^2 \, - \, 1 ) u_\lambda 
     \; = \; 0
      \quad \mathrm{in} \quad Q(T),
\label{EQ:GLHF-Org}
\\[2mm]
u_\lambda
 \; = \; u_0
  \quad \mathrm{on} \quad \partial Q(T).
\label{EQ:GLHF-Bdry}
\end{eqnarray}
If you notice that the nonlinear term of 
 $\lambda^{1-\kappa}$
  $\dot \chi \bigl( \bigl( | u_\lambda |^2 \, - \, 1 \bigr)^2 \bigr)$
   $\bigl( | u_\lambda |^2 \, - \, 1 \bigr) u_\lambda$
    is bounded, 
     Banach's fixed point theorem can state
      the unique existence of the mapping $u_\lambda$ on $Q(T)$ with
\begin{enumerate}
\item[(a)]
$u_\lambda \, \in \, C^\infty (Q(T))$,
\item[(b)]
\text{\eqref{EQ:GLHF} is fulfilled in } $Q(T)$,
\item[(c)]
$u_\lambda (t,x) \, - \, u_0 (x)$
 $\in$
  $\tc_B (\Omega;\mathbb{R}^{D+1})$ 
$\;\, \text{for almost every}$ 
 $\, t \;\mathrm{in}\; (0,T)$,
\item[(d)]
$
\lim_{t \searrow 0}
 || u_\lambda (t, \cdot) \, - \, u_0 (\cdot) ||_{L_B^2 (\Omega)}
  \; = \; 0$.
\end{enumerate}
\label{DEF:GLHF}
\par
An observation by K.Horihata~\cite[Theorem 2.4]{horihata}
 and an application of a maximal principle on an unbounded domain
  which referred to F.John~\cite[Chap 7]{john}
   prevails $| u_\lambda |$ $\le$ $1$ in $Q(T)$.
Thus we rewrite \eqref{EQ:GLHF-Org} to
\begin{equation}                                                                                                                                                                  
\dfrac{\partial u_\lambda}{\partial t} \, - \, 
 \triangle u_\lambda
  \, + \, \lambda^{1-\kappa}
   \bigl( | u_\lambda |^2 \, - \, 1 \bigr) u_\lambda                                                                                                                              
    \; = \; 0
     \quad \mathrm{in} \quad Q(T).                                                                                                                                                
\label{EQ:GLHF}
\end{equation}
\par
The benefits of my scheme is to easily derive
\begin{equation}
\lint_{Q(T)} \lambda^{1-\kappa}
 ( | u_\lambda |^2 \, - \, 1 )^2 \, dz
\; =  \; O ( 1/\log\lambda)
 \quad \mathrm{as} \; \lambda \nearrow \infty.
\notag
\end{equation}
\par
We outlook the history on a boundary regularity of harmonic mapping or harmonic heat flows:
For the energy minimizing map, R.Schoen and K.Uhlenbeck~\cite{schoen-uhlenbeck} has established 
 a boundary regularity for an appropriate boundary condition.
  A similar results on minimizers of a certain functional
   can be found in J.Jost and M.Meier~\cite{jost-meier}.
J.Qing~\cite{qing} proved the boundary regularity of weakly harmonic maps from surfaces with the boundary.
 C.Poon~\cite{poon} constructed a smooth harmonic map $u_a$ between 
  $B^3$ and $\mathbb{S}^2$ except a prescribed point $a$ $\in$ $\overline{B}^3$,
   noncontinuous at $a$ and $u_a$ $=$ $x$ on $\partial B^3$.
\par
On the other hand, the papers discussed boundary regularity results on
 a harmonic heat flow are not so many.
We refer it to Y.Chen \cite{chen-91} or Y.Chen-F.-H.Lin \cite{chen-lin} 
 or Y.Chen-J.Li-F.H.Lin \cite{chen-li-lin}
  or C.Y.Wang \cite{wang-99}.
Since the maps $u_a$ by C.Poon \cite{poon} is 
 also a non-smooth weakly harmonic heat flow,
  it may be reasonable to discuss 
   a partial regularity for our harmonic heat flow
    and a regularity under a certain imposition.
\par
This paper has two folds:
 The first constructs a few energy inequalities
  on $\overline{\Omega}$ 
   and next discusses a partial regularity 
    for our harmonic heat flow near the boundary;
The proof proceeds as in the one of the author's former paper
  of K.Horihata~\cite{horihata}
   combined with the one by Y.Chen~\cite{chen-91}.
Thereafter by utilizing the reverse Poincar\'e inequality,
 we will prove that 
  the WHHF is actually smooth except on a small set
   called \lq\lq \bf{singular set}.\rq\rq \rm \enspace
    More precisely we assert
\begin{Thm}{\rm{(Partial Regularity Theorem).}}\label{THM:Main-1}
Let $d$ be a positive integer larger than $2$.
For a mapping $u_0$ $\in$ $H_B^{1,2}$ $( \Omega \, ; \mathbb{S}^D )$
  $\cap$
   $(H_B^{2-1/q,q} \cap H_B^{[(d+1)/2]+1,p})$ $( \Omega\setminus(\complement\Omega)_{\delta_0} ; \mathbb{S}^D)$
    with a positive number sufficiently small $\delta_0$,
     and positive numbers $q$ and $p$ respectively greater than or equal to $d+2$
      and greater than $2$,
there exists a WHHF
 and it is smooth on a certain relative open set in $\overline{Q}(T)$ off a set called {\bf{sing}}.
  The set has the finite $(d-\gamma_0)$-dimensional Hausdorff measure
   with respect to the parabolic metric,
where $\gamma_0$ is a small positive number depending only on 
 $u_0$, $d$ and $Q(T)$.
The WHHF also holds
\allowdisplaybreaks\begin{align}
&
\lint_{t_0 - (2R_1)^2}^{t_0 - R_1^2}
 \lint_\Omega 
  | \nabla u |^2 G_{z_0} \, dz
\, + \, 2 \lint_{R_1}^{R_2} \, dR
 \lint_{t_0 - (2R)^2}^{t_0 - R^2}
  \, dt
   \lint_\Omega
\Bigl\vert
 \frac {\partial u_\lambda} {\partial t}
  \, - \, \frac{x-x_0}{2\sqrt{t_0-t}}
   \cdot \nabla u 
    \Bigr\vert^2 
G_{z_0} \, dx
\notag
\\
&
\; \le \; C 
 (R_2^{\mu_0} \, - \, R_1^{\mu_0})
  \lint_{t_0 - (2R_2)^2}^{t_0 - R_2^2} \, dt
   \lint_\Omega
    | \nabla u |^2 G_{z_0} \, dx
\label{INEQ:Boundary-Monotonical-1}
\, + \, 
 C (R_2 \, - \, R_1)
\end{align}
for any positive numbers $R_1$ and $R_2$ with $R_1 \le R_2$ 
 and an arbitrary point $z_0$ $=$ $(t_0,x_0)$ in $\overline{Q}(T)$
  satisfying $R_2$ $<$ $\sqrt{t_0/4}$,
where two positive constants $\mu_0$ and $C(\mu_0)$
 has the relation of $C(\mu_0)$ $\nearrow$ $\infty$
  as $\mu_0$ $\searrow$ $0$
and in addition
\allowdisplaybreaks\begin{align}
&
\lint_{P_R (z_0)} | \nabla u (z) |^2 \, dz
 \; \le \; \frac {C}{R^2} \lint_{P_{2R} (z_0)}
  \vert u (z) \, - \, h_0 (x) \vert^2 \, dz
\notag
\\
&
\, + \, \lint_{P_R (z_0)} 
 ( \vert \nabla^{[(d+1)/2]+1} h_0 (x) \vert^2 
  \, + \, \vert \nabla h_0 (x) \vert^2 )
   \, dz
\end{align}
for any parabolic cylinder $P_{2R}$ $(z_0)$,
 where a function $G_{z_0}$ means a backward heat kernel indicated by
\begin{equation}
G_{z_0} (t,x) \; = \; \frac 1{\sqrt{4\pi (t_0-t)}^d}
 \exp \Bigl( - \frac {|x-x_0|^2}{4(t_0-t)} \Bigr)
  \quad ((t,x) \, \in \, (0,t_0) \times \mathbb{R}^d).
\notag
\end{equation}
\end{Thm}
\par
The subsequent theorems will take up
 two sufficient conditions on the whole domain's regularity
  for our {\it{WHHF}}.
The one is the boundary energy smallness at a large time and
 the another is a target restriction
  on the initial mapping called
   \lq\lq one-sided condition\rq\rq.
    We refer it to M.Giaquinta~\cite[p.237, Theorem 3.2]{giaquinta}.
\begin{Thm}{\rm{(First Regularity Theorem).}}\label{THM:Main-2}
Assume that our domain is bounded with \eqref{Cond:B}
and there exists some positive number $\epsilon_0$ such that
 for all point $z_0$ $=$ $(t_0,x_0)$ and some positive number $R_0$
  less than $\sqrt{t_0}/2$,
   if we have
\allowdisplaybreaks\begin{align}
&
\frac {\exp^{(4R_0)^{\mu_0}}}{R_0^2}
 \Bigl[
  \frac {\exp^{-4(d-2)/d_0^2}}{t_0^{(d-2)/2}}
   \lint_{\Omega} | \nabla u_0 |^2 \, dx
\Bigr.
\label{INEQ:Energy-Decay-Bdry}
\\
&
\Bigl.
\, + \,
 \lint_0^{t_0-R_0^2} \, dt
  \lint_{\partial \Omega} | \nabla_\tau u_0 |^2 G_{z_0}
   \bigl(d_{x_0} \, + \, \frac {4(t_0-t)}{d_0^2} \bigr)
    \, d\mathcal{H}_x^{d-1}
     \Bigr]
\, + \, C(\mu_0) R_0
\; \le \; \epsilon_0^2,
\notag
\end{align}
where $\mu_0$ and $C(\mu_0)$ are mutual relevant positive constants
 satisfying $C(\mu_0)$ $\searrow$ $\infty$ as $\mu$ $\nearrow$ $0$.
\begin{equation}
d_{x_0} (x) \; = \; 1 \, + \, \frac {|x-x_0|^2}{d_0^2}
 \quad \mathrm{and}
  \quad d_0 \; = \; \diam (\Omega).
\notag
\end{equation}
Then our WHHFs are smooth on a neighbourhood of $z_0$ $=$ $(t_0,x_0)$.
\end{Thm}
\begin{Thm}{\rm{(Second Regularity Theorem).}}\label{THM:Main-3}
If the range of $u_0${\rm{:}} $u_0$ $(\overline\Omega)$ is compactly contained in $\mathbb{S}_+^D$
 aftermath of a suitable rotation if necessary,
  then the range of our heat flow $u(\overline{Q}(T))$ is done there 
   and they are smooth in $(0,T)\times\overline{\Omega}$.
\end{Thm}
\vskip 9pt
\par
We list a glossary of notation below{\rm{:}}
 We adopt it from K.Horihata~\cite{horihata} and 
  add it up the new symbols especially related to the boundary.
The symbols of vector and the definition on
 the function class can be seen therein.
\vskip 9pt
\begin{center}
\underline{\textit{Notation}}
\end{center}
\vskip 9pt
\renewcommand{\labelenumi}{(\roman{enumi})}
\begin{enumerate}
\item
$\mathbb{R}_\pm^d$ $=$
 $\{x = (x^\prime,x_d) \in \, \mathbb{R}^d
  \; ; \;
   x_d \, \gtrless \, 0 \}$
    and 

     $\mathbb{R}_0^d$ $=$
      $\{x = (x^\prime,x_d) \in \, \mathbb{R}^d
       \; ; \; x_d \, = \, 0 \}$
        $(\backsimeq \mathbb{R}^{d-1})$.
\item
For any set $A$ $\subset$ $\mathbb{R}^d$ and 
 any positive number $\delta$,
  $A_\delta$ $=$ $\{ x \in \mathbb{R}^d \, ; \, 
    \inf_{y \in \Omega} |x-y| \, < \, \delta \}$.
\item
The symbol $\Omega$ is a domain in $\mathbb{R}^d$                                           
 and $\; \partial \Omega$ the boundary of $\Omega$.
  Moreover for a set $A$ $\subset$$\mathbb{R}^d$,
   the symbol $A_\pm$ respectively means 
    $A \cap \Omega$ and $A \cap \complement\Omega$
     and $A_0$ does $A \cap \partial \Omega$.
\item
A vector $\nu$ and $\tau$ respectively denotes
 the unit outer normal and the tangential field along $\partial \Omega$.
\item
${Q(T)}=(0,T) \times \Omega$.
 $\; \partial {Q(T)} \, = \,$
  $[0,T) \times \partial \Omega$
   $\cup$
    $\{0\} \times \Omega$.
\item
$\mathbb{S}^D \, = \,$ 
 $\{y = (y^1,y^2,\ldots, y^{D+1}) \, \in \, \mathbb{R}^{D+1};$ 
  $|y| \, = \, \sqrt{\sum_{i=1}^{D+1} (y^i)^2} \, = \, 1\}$.
\item
$\mathbb{S}_\pm^D \, = \,$
 $\{y = (y^\prime,y^{D+1}) = (y^1,y^2,\ldots, y^{D+1}) \, \in \, \mathbb{S}^D$;
  $y^{D+1} \, \gtrless \, 0 \}$.
\item
$\mathbb{S}_0^{D} \, = \,$
 $\{y = (y^\prime,y^{D+1}) = (y^1,y^2,\ldots, y^{D+1}) \, \in \, \mathbb{S}^D$;
  $y^{D+1} \, = \, 0 \}$

   $(\backsimeq \partial B_1^D (0))$.
\item
Set points $x$ $=$ $(x_1,x_2,\ldots,x_d)$,
 $x_0$ $=$ $(x_{0,1},x_{0,2},\ldots,x_{0,d})$ and
  $z_0 \, = \, ( t_0 , x_0) \in \mathbb{R}^{d+1}$;
   Then indicate a ball $B_r (x_0)$, a parabolic cylinder $P_r (z_0)$ by
\allowdisplaybreaks\begin{align*}
B_r (x_0) &
 \, = \, \{ x \in \mathbb{R}^d \, ; \,  | x \, - \, x_0 | < r \}, \\
P_r (z_0) & \, = \, 
 \{ z=(t,x) \in Q(T)  \, ; \, t_0 - r^2 < t < t_0 + r^2,
  | x \, - \, x_0 | < r\}.
\end{align*}
In $B_r (x_0)$ and $P_r (z_0)$,
 the points of $x_0$ and $z_0$ will be often abbreviated
  when no confusion may arise.
\item
Letter $C$ denotes a generic constant.
 By the letter $C(B)$, it means that a constant depends only on a parameter $B$.
\item
$\kappa (t)$ is $\arctan (t) / \pi$.
\item
For any point $z_0$ $=$ $(t_0,x_0)$ in $\mathbb{R}^{d+1}$,
 the back-ward heat kernel $G_{z_0}$ designates
\begin{equation}
G_{z_0} (t,x) \; = \; \frac 1{\sqrt{4\pi (t_0-t)}^d}
 \exp \Bigl( - \frac {|x-x_0|^2}{4(t_0-t)} \Bigr)
 \quad ((t,x) \, \in \, (0,t_0) \times \mathbb{R}^d).
\notag                                                                                                                                                                            \end{equation}
\item
For the nabla:
\begin{math}
\nabla \; = \; \bigl(
 \partial/\partial y_1,
  \partial/\partial y_2,
   \ldots,
    \partial/\partial y_d
     \bigr),
\end{math}
the differential operators $\nabla_\nu$ and $\nabla_\tau$ denote
\allowdisplaybreaks\begin{align*}
&
\nabla_\nu \; = \; \nu \cdot \nabla, \;
\nabla_\tau \; = \; \nabla \, - \, \nu \nu \cdot \nabla.
\notag
\end{align*}
\end{enumerate}
\par
A forthcoming paper extends the results here
 and the ones of K.Horihata~\cite{horihata} to a heat flow between
  any Riemannian manifold under a certain assumption.

%
%
\setcounter{chapternumber}{2}\setcounter{equation}{0}
\renewcommand{\theequation}%
             {\thechapternumber.\arabic{equation}}
\section{\enspace GLHF.}
\label{P:GLHF}
The chapter introduces three fundamental energy inequalities:
 The first is a monotonical inequality
  and the second is a local decay energy estimate,
   which is can be regarded as a variant of a monotonical inequality,
    and the final is a hybrid inequality.
\subsection{\enspace Two Monotonical Inequalities.}
We introduce former two inequalities. 
 The one is it by Y.Chen and F.H.Lin~\cite{chen-lin} and
   the another is a hybrid type inequality 
    combined with the energy decay estimate
     referred to K.Horihata~\cite[Theorem 2.6]{horihata} and
      the usual monotonical inequality:
\begin{Thm}{\rm{(Monotonical Inequalities).}}\label{THM:Mon}
There exist constants $C$ and $C(\mu_0)$
 with $C(\mu_0)$ $\nearrow$ $\infty$ as $\mu_0$ $\searrow$ $0$
  such that the following holds
   for any point $z_0 \, = \, (t_0,x_0)$ $\, \in \, Q (T)$ and 
    positive numbers $R$, $R_1$ and $R_2$ with $0 < R_1 < R_2$
     with $\max(R_2,R)$ $\, \le \, \sqrt{t_0}/2$,
\allowdisplaybreaks\begin{align}
&
\lint_{t_0 - 4R_1^2}^{t_0 - R_1^2} \, dt \lint_{\Omega} 
 \mathbf{e}_\lambda \, G_{z_0} \, dx      
\, + \, \lint_{R_1}^{R_2} \, dt
 \lint_{\Omega}
  \left| \frac {\partial u_\lambda}{\partial t}
   \; - \;
    \frac {x \, - \, x_0}{2(t \, - \, t_0)} \cdot \nabla u_\lambda \right|^2
     G_{z_0} \, dx
\notag
\\
&
\; \le \; C \exp{(R_2^{\mu_0} \, - \, R_1^{\mu_0})}
 \lint_{t_0 - 4R_2^2}^{t_0 - R_2^2} \, dt \lint_{\Omega}
  \mathbf{e}_\lambda \, G_{z_0} \, dx
\notag
\\
&
\, + \, C(\mu_0) (R_2 \, - \, R_1)
\label{INEQ:Mon}
\end{align}
and in addition under the hypothesis \eqref{Cond:B},
%
%
\allowdisplaybreaks\begin{align}
& R^2 \lint_{\{t_0 - R^2\} \times \Omega}
 \mathbf{e}_\lambda 
  G_{(t_0-R^2,x_0)} \, dx
\; \le \; 
 \frac {e^{-4(d-2)(t_0 - R^2)/d_0^2}}{(2\sqrt{t_0})^d}
  \lint_{\Omega} | \nabla u_0 |^2 \, dx
\notag
\\
&
\, + \, C
 \lint_0^{t_0-R^2} \, dt
  \lint_{\partial \Omega} | \nabla_\tau u_0 |^2
   \, G_{z_0} \, d\mathcal{H}_x^{d-1}
\label{INEQ:Energy-Decay-Bdry-GLHF}
\end{align}
with 
\allowdisplaybreaks\begin{align}
&  
G_{z_0} (t,x) \; = \; \frac 1{\sqrt{4\pi (t_0-t)}^d}
 \exp \Bigl( - \frac {|x-x_0|^2}{4(t_0-t)}\Bigr)
  \quad ((t,x) \, \in \, (0,t_0) \times \mathbb{R}^d),
\notag
\\
&
d_{x_0} (x) \; = \; 1 \, + \, \frac {|x-x_0|^2}{d_0^2}
 \quad \mathrm{with}
  \quad d_0 \; = \; \diam (\Omega).
\notag
\end{align}
\end{Thm}
\vskip 9pt
\noindent{\underbar{{Proof of Theorem \ref{THM:Mon}}}.}
\rm\enspace
\vskip 6pt
Set 
\allowdisplaybreaks\begin{align}
&   
\nabla_{|x-x_0|} \; = \; \frac {x-x_0}{|x-x_0|} \cdot \nabla, \;
\nabla_{\tau_0} \; = \; 
 \nabla \, - \, \frac {x-x_0}{|x-x_0|} \nabla_{|x-x_0|}.
\end{align}
Since \eqref{INEQ:Mon} is well-known, 
 we only show \eqref{INEQ:Energy-Decay-Bdry-GLHF}:
Multiplying \eqref{EQ:GLHF} by
$$
2 (t_0-t) \frac {\partial u_\lambda}{\partial t}
 G_{z_0} d_{x_0}
  \quad \mathrm{and} \quad
   (x-x_0) \cdot \nabla u_\lambda G_{z_0}
    ( d_{x_0} \, + \, 4(t_0-t)/d_0^2)
$$
and integrating it over $\Omega$, we obtain
\allowdisplaybreaks\begin{align}
&
2 \lint_{\Omega}
 \left| \frac {\partial u_\lambda}{\partial t}
  \right|^2 (t_0-t) G_{z_0} d_{x_0} \, dx
\, + \, 2
 \frac d {dt}
  \lint_{\Omega} \mathbf{e}_\lambda
   (t_0 - t) \, G_{z_0} d_{x_0} \, dx
\notag
\\
&
\, - \, (d-2) 
 \lint_{\Omega} \mathbf{e}_\lambda
  \, G_{z_0} d_{x_0} \, dx
\, + \, 
 \lint_{\Omega} \mathbf{e}_\lambda 
  \frac {|x-x_0|^2}{2(t_0 - t)} \, G_{z_0} d_{x_0} \, dx
\notag
\\
&
\, + \, \lambda^{1-\kappa}
 \lint_{\Omega}
  ( |  u_\lambda |^2 \, - \, 1 )^2
   \, \frac {|x-x_0|^2}{8(t_0 - t)} \, G_{z_0} \, d_{x_0} \, dx
\notag
\\
&
\, - \, \lint_{\Omega}
 \Bigl\langle \frac {\partial u_\lambda}{\partial t},
  (x-x_0) \cdot \nabla u_\lambda
   \Bigr\rangle 
\, G_{z_0} 
 \bigl( d_{x_0} \, - \, \frac {4(t_0-t)}{d_0^2}
  \bigr) \, dx
\label{INEQ:ProofOfMono-1}
\\
&
\, + \, \dot\kappa \log \lambda \, \frac {\lambda^{1-\kappa}}2
 \lint_{\Omega}
  ( |  u_\lambda |^2 \, - \, 1 )^2
   \, (t_0 - t) \, G_{z_0} \, d_{x_0} \, dx
    \; = \; 0,
\notag
\\
&
\, - \, \lint_{\Omega}
 \Bigl\langle \frac {\partial u_\lambda}{\partial t},
  (x-x_0) \cdot \nabla u_\lambda,
   \Bigr\rangle
\, G_{z_0} \, 
 \Bigl( d_{x_0} \, + \, \frac {4(t_0-t)}{d_0^2} \Bigr)
  \, dx
\notag
\\
&
\, + \, (d-2)
 \lint_{\Omega} \mathbf{e}_\lambda
  \, G_{z_0} 
   \Bigl( d_{x_0} \, + \, \frac {4(t_0-t)}{d_0^2} \Bigr)
    \, dx
\notag
\\
&
\, + \,
 \lint_{\Omega}
  ( | \nabla_{|x-x_0|} u_\lambda |^2 \, - \, | \nabla_{\tau_0} u_0 |^2 )
   \, \frac {|x-x_0|^2}{4(t_0-t)}
    \, G_{z_0} \, d_{x_0} \, dx      
\label{INEQ:ProofOfMono-2}
\\
&
\, + \, \frac {\lambda^{1-\kappa}}2
 \lint_{\Omega}
  ( |  u_\lambda |^2 \, - \, 1 )^2
   \, \frac {|x-x_0|^2}{t_0-t}
    G_{z_0} \, 
     \bigl( d_{x_0} \, + \, \frac {4(t_0-t)}{d_0^2} \bigr)
      \, dx
\notag
\\
&
\, - \, \frac {\lambda^{1-\kappa}}8
 \lint_{\Omega}
  ( |  u_\lambda |^2 \, - \, 1 )^2
   \, \frac {|x-x_0|^2}{t_0-t}
    G_{z_0} \, d_{x_0} \, dx
\notag
\\
&
\, + \, \frac 12 \lint_{\partial\Omega}
 ( \langle x-x_0,\nu \rangle 
  \vert \nabla_\nu u_\lambda \vert^2
   \, + \, 2 
    \langle \nabla_\nu u_\lambda,
     (x-x_0) \cdot \nabla_\tau u_0
      \rangle
       \, - \, \langle x-x_0,\nu \rangle
        \vert \nabla_\tau u_0 \vert^2
         )      
\notag
\\
&
\quad \quad \times G_{z_0} \,
 \Bigl( d_{x_0} \, + \, \frac {4(t_0-t)}{d_0^2} \Bigr)
  \, d\mathcal{H}_x^{d-1}
\; = \; 0.
\notag
\end{align}
\par
Summarizing \eqref{INEQ:ProofOfMono-1}, \eqref{INEQ:ProofOfMono-2}
 and using \eqref{Cond:B}, we arrive at
\allowdisplaybreaks\begin{align}
&
\frac d {dt}
 \lint_{\Omega} \mathbf{e}_\lambda
  \, (t_0-t) d_{x_0} G_{z_0} \, dx
\, + \, \frac{4(d-2)}{d_0^2} \lint_{\Omega} \mathbf{e}_\lambda
 \, (t_0-t) G_{z_0} d_{x_0} \, dx
\notag
\\
&
\, + \, 2 \lint_{\Omega}
 \Bigl| \frac {\partial u_\lambda}{\partial t} \sqrt{t_0-t}
  \, - \, \frac {x-x_0}{2\sqrt{t_0-t}} \cdot \nabla u_\lambda
   \Bigr|^ 2G_{z_0} d_{x_0} \, dx
\, + \, 
 \frac{\theta_0}2
  \lint_\Omega
   \vert x^\prime \, - \, x_0^\prime \vert^2
    \vert \nabla_\nu u_\lambda \vert^2
     G_{z_0} \, d_{x_0} \, dx
\notag
\\
&
\; \le \; C
 \lint_{\partial \Omega}
  \vert \nabla_\tau u_0 \vert^2
G_{z_0} \, 
 \Bigl(
  d_{x_0} \, + \, \frac {4(t_0-t)}{d_0^2} 
   \Bigr)
\, d\mathcal{H}_x^{d-1}.
\label{INEQ:ProofOfMono-3}
\end{align}
\par
A multiplier of \eqref{INEQ:ProofOfMono-3} by $e^{4(d-2)t/d_0^2}$
 and an integration on it from $0$ to $t_0 - R^2$
  imply our claim.
%
%
\subsection{\enspace Hybrid type Inequality for GLHF.}
We demonstrate the reverse Poincar\'{e} inequality.
 To this end we prepare the below:
  Let the mapping $h_0$ be the solution to
\begin{equation}
\left\{
\begin{array}{rl}
- \triangle h_0 & \; = \; 0 
\quad \mathrm{in} \quad \Omega, \\
   h_0 & \; = \; u_0
\quad \mathrm{on} \quad \partial\Omega.
\label{EQ:HF}
\end{array}
\right.
\end{equation}
Then we claim
\begin{Thm}{\rm{(Hybrid Inequality).}}\label{THM:HI}
For any positive number $\epsilon_0$ and
 any point $z_0$ in $Q(T)$, 
  there exists a positive constant
   $C(\epsilon_0)$ satisfying $C(\epsilon_0)$
    $\nearrow$ $\infty$ as $\epsilon_0 \searrow 0$
     such that the inequality
\allowdisplaybreaks\begin{align}
& \lint_{P_{R} (z_0) \cap Q(T)} \mathbf{e}_\lambda (z) \, dz
 \, \le \; \epsilon_0 \lint_{P_{2R} (z_0) \cap Q(T)} \mathbf{e}_\lambda (z) \, dz
\, + \, 
 \frac {C(\epsilon_0)}{R^2} \lint_{P_{2R} (z_0) \cap Q(T)}
  \vert u_\lambda (z) \, - \, h_0 (x) \vert^2 \, dz
\notag
\\
&
\, + \, \lint_{P_{2R} (z_0) \cap Q(T)}
 ( \vert \nabla^{[(d+1)/2]+1} h_0 (x) \vert^2 \, dz
  \, + \, 
   \vert \nabla h_0 (x) \vert^2 ) \, dz
    \, + \, o(1)
     \quad ( \lambda \nearrow \infty)
\label{INEQ:HI}
\end{align}
holds for any parabolic cylinder $P_{2R} (z_0)$.
\end{Thm}
\par
Preliminary we rewrite \eqref{EQ:GLHF} to it by the flatten-out coordinate:
 Set $v_\lambda (t,y)$ on $\Phi_{x_0}$ $(U (x_0) \cap \Omega)$
  by $v_\lambda (t,y)$ $=$ $u_\lambda$ $(t, \Phi_{x_0}^{-1}(y))$.
   Then \eqref{EQ:GLHF} becomes
\begin{eqnarray}
\dfrac{\partial v_\lambda}{\partial t} \, - \,
 L(v_\lambda)
  \, + \, \lambda^{1-\kappa}
   \bigl( | v_\lambda |^2 \, - \, 1 \bigr) v_\lambda 
    \; = \; 0
     \; \mathrm{in} \; (0,T) \times \Phi_{x_0} (U (x_0) \cap \Omega)
\label{EQ:GLHF-Flat}
\\[2mm]
\notag
\end{eqnarray}
with
\allowdisplaybreaks\begin{align}
&
D_i \; = \; 
 \frac {\partial}{\partial y_i}
  \, - \, a_i \frac {\partial}{\partial y_d},
   \quad a_i \; = \; \frac {\partial \phi_{x_0}}{\partial y_i}
    \quad ( i \, = \, 1,2,\ldots,d-1),
\notag
\\
&
D_d \; = \; \frac {\partial}{\partial y_d}, \;
L \; = \; \sum_{i=1} ^{d-1} D_i^2 \, + \, D_d^2,
\notag
\\
&
D_\nu \; = \; \sum_{i=1}^d \frac {y^i}{|y|} D_i,
 \;
  D_\tau \; = \; D \, - \, \frac {y}{|y|} \, D_\nu,
\notag
\\
&
\triangle_\tau \; = \; \triangle \, - \, 
 \frac 1{\rho^{d-1}} \frac {\partial}{\partial \rho}
  \Bigl( \rho^{d-1} \frac{\partial}{\partial \rho} \Bigr).
\end{align}
Choose $B_R (0)$ $\subset$ $\Phi_{x_0} (U (x_0))$
 and induce the mapping $h_{x_0}$ by 
\begin{equation}
h_{x_0} (t,y) \; = \; h_0 (\Phi_{x_0}^{-1} (y)).
\end{equation}
\par
We extend the mapping $v_\lambda \,- \, h_{x_0}$ in 
 $B_R (0) \cap \mathbb{R}_+^d$ to the one in $B_R (0)$ by
\begin{equation}                                                                                                                                                                  
\left\{                                                                                                                                                                           
\begin{array}{rlc}                                                                                                                                                             
( v_\lambda \, - \, h_{x_0} ) (t,(y^\prime,y_d)) &
 \quad \mathrm{if} \quad y \, \in \, B_R (0) \cap \mathbb{R}_+^d,
\notag
\\[1mm]
0 &
\quad \mathrm{if} \quad y \, \in \, B_R (0) \cap \mathbb{R}_0^d,
\notag
\\[1mm]
- ( v_\lambda \, - \, h_{x_0} ) (t,(y^\prime,-y_d)) &
\quad \mathrm{if} \quad y \, \in \, B_R (0) \cap \mathbb{R}_-^d
\notag
\end{array}
\right.                                                                                                                                                                           
\label{EQ:GLHF-Flat}
\end{equation}           
and denote it by the same symbol $v_\lambda \, - \, h_{x_0}$:
 Note that the mapping $v_\lambda \, - \, h_{x_0}$
  belongs to $(C^1 \cap W^{2,2})$ $(B_R (0) \, ; \, \mathbb{R}^{D+1})$
   and it satisfies the following identity:
\begin{eqnarray}
\lint_{B_R(0)}
 \Bigl\langle 
  \dfrac{\partial v_\lambda}{\partial t},
   \phi 
    \Bigr\rangle \, dx
\, + \, \lint_{B_R(0)}
 \langle D (v_\lambda - h_{x_0} ), D \phi \rangle \, dx
\notag
\\
\, + \, \lambda^{1-\kappa} \lint_{B_R(0)}
 \bigl( | v_\lambda |^2 \, - \, 1 \bigr) 
  \langle v_\lambda, \phi \rangle \, dx
   \; = \; 0
\label{EQ:GLHF-Boundary}
\end{eqnarray}
for any mapping $\phi$ $\in$ $C_0^\infty$ 
 $( B_R (0) \, ; \, \mathbb{R}^{D+1})$
  with
\allowdisplaybreaks\begin{align}
&     
\langle D (v_\lambda - h_{x_0} ), D \phi \rangle 
 \; = \; 
  \sum_{i=1}^d \sum_{j=1}^{D+1}
   D_i (v_\lambda - h_{x_0} )^j D_i \phi^j
\notag
\\
&
\Bigl\langle \frac {\partial v_\lambda}{\partial t},
 \phi \Bigr\rangle
  \; = \; \sum_{j=1}^{D+1}
   \frac {\partial v_\lambda^j}{\partial t} \phi^j,
\;
 \langle v_\lambda, \phi \rangle
  \; = \; \sum_{j=1}^{D+1}
   v_\lambda^j \phi^j
\notag
\end{align}
to verify
\begin{Lem}\label{LEM:GLQE}
For any balls $B_{\rho_1}$ and $B_{\rho_2}$ 
 with $B_{\rho_1}$ $\, \subset \,$ $B_{\rho_2}$
  $\, \subset \,$ $B_R (0)$
   and any cylinders $P_{\rho_1}$ and $P_{\rho_2}$
    with $P_{\rho_1}$ $\, \subset \,$ $P_{\rho_2}$
     $\, \subset \,$ $P_R (0)$, the following holds{\rm{:}}
\allowdisplaybreaks\begin{align}
&
\lint_{P_{\rho_1}} \lambda^{1-\kappa}
 ( 1 \, - \, | v_\lambda |^2 ) \, dz
\; \le \; C \lint_{P_{\rho_2}}
 \mathbf{e}_\lambda \, dz
\notag
\\
&
\, + \, \frac C{(\rho_2 - \rho_1)^2}
 \lint_{P_{\rho_2}}
  ( 1 \, - \, | v_\lambda |^2 ) \, dz,
\label{INEQ:GLQE-1}
\\
&
\lambda^{1-\kappa}
 \lint_{B_{\rho_1}}
  ( 1 \, - \, | v_\lambda |^2 ) \, dy
\; \le \; C \lint_{B_{\rho_2}}
 \mathbf{e}_\lambda \, dy
\notag
\\
&
\, + \, \frac {C}{\rho_2 - \rho_1}
 \lint_{B_{\rho_2}}
  \vert D v_\lambda \vert \, dy
\, + \, C \lint_{B_{\rho_2}}
 \Bigl\vert \frac {\partial v_\lambda}{\partial t }\Bigr\vert
  \, dy.
\label{INEQ:GLQE-2}
\end{align}
\end{Lem}
\par
Secondly we list symbols and auxiliary mapping employed only here.
 Give $L_\lambda$ by $[ \log (\lambda / h(\lambda)) / \log 2 ] + 1$,
  where a function $h(\lambda)$ $( \lambda \in \mathbb{R})$ is positive
   satisfying $h(\lambda)$ $\searrow$ $0$ $( \lambda \nearrow \infty)$.
    We then  introduce the decompositional convention: 
Let $r$ be a positive number less than $\epsilon_0^2 R$,
 and put
\allowdisplaybreaks\begin{align}
&
\triangle \rho_l \; = \; C_1 \epsilon_0^4 r (1/2)^l
 \quad ( l = 1,2, \ldots, L_\lambda ),
\notag
\\
\rho_l \; = \; & 
\left\{
\begin{array}{ll}
0 & \quad ( l = 0 )
\\[2pt]
( 1 - \epsilon_0^4 ) r
 \, + \, 
  C_1 \epsilon_0^4 r \sum_{j=1}^{l} (1/2)^j
& \quad ( l = 1, 2, \ldots, L_\lambda )
\end{array}
\right.
\notag
\\
&
\textrm{with} \;
 C_1 \, = \, 
  \Bigl(\sum_{l=1}^{L_\lambda} (1/2)^l \Bigr)^{-1}.
\notag
\end{align}
Fix $t$ $\in$ $(-r^2,r^2)$ and
 choose numbers $r_l$ $( l \, = \, 1,2, \ldots, L_{\lambda} -1 )$
  and $\triangle r_l$
   $( l \, = \, 1,2, \ldots, L_{\lambda} )$ 
    so that they are given by
\allowdisplaybreaks\begin{align}
&
\frac{\triangle \rho_l}{3}
 \lint_{\{r_l\} \times \mathbb{S}^{d-1}}
  | D_\nu ( v_\lambda \, - \, h_{x_0} ) ( t,y ) |^2 \, %
   d\mathcal{H}_y^{d-1}
\; = \; 
 \lint_{\rho_{l} - \triangle \rho_l/6}
      ^{\rho_{l} + \triangle \rho_l/6}
   \rho^{d-1} \, d\rho
    \lint_{\mathbb{S}^{d-1}}
     |  D_\nu ( v_\lambda \, - \, h_{x_0} ) ( t,y ) |^2 \, d\omega_{d-1},
\notag
\\
&
r_0 \, = \, 0, \;
 r_{L_\lambda} \, = \, \rho_{L_\lambda}
  \; = \; r,
\triangle r_l \; = \; r_l \, - \, r_{l-1}.
\notag
\end{align}
In addition a positive number $\hat{r}_l$
 $( l \, = \, 0, 1, \ldots, L_\lambda )$ denotes
\allowdisplaybreaks\begin{align}
&
\widehat{r}_l \; = \; 
\left\{
\begin{array}{lcl}
0 & & \quad ( l = 0 )
\\[2pt]
\dfrac {{r}_{l} \, + \, {r}_{l+1}}2
& & \quad ( l \, = \, 1, 2, \ldots, L_\lambda -1 )
\\[2pt]
r & & \quad ( l = L_\lambda ).
\end{array}
\right.
\notag
\end{align}
\par
Next introduce a mapping $f_{\lambda}$ which is the solution to
\allowdisplaybreaks\begin{align}
&
\left\{
\begin{array}{rcl}
D_\nu f_{\lambda}
 \, + \, r \triangle_\tau f_{\lambda} \; & = \; & 0
  \quad \mathrm{in} \quad ( 0, r ) \times \mathbb{S}^{d-1},
\\[2mm]
f_{\lambda} \; & = \; & v_{\lambda} \, - \, h_{x_0}
 \quad \mathrm{on} \quad \{ {r} \} \times \mathbb{S}^{d-1}.
\end{array}                                                                                                                                                                      
\right.                                                                                                                                                                          
\label{EQ:Support}                                                                                                                                                            
\end{align}
\par
Designate five sorts of annuls:
\allowdisplaybreaks\begin{align}
&
\widehat{T}_l \; = \;
 [ \widehat{r}_{l-1}, \widehat{r}_{l}) \, \times \, {\mathbb{S}^{d-1}},
  \quad
T_l \; = \; [ \, r_{l-1}, r_{l}) \, \times \, {\mathbb{S}^{d-1}},
\notag
\\
&
T_l^1 \; = \; [r_{l-1}, r_{l-1} \, + \, \triangle r_l/3)
 \, \times \, \mathbb{S}^{d-1},                                                                                                                                                  
\notag
\\
&
T_l^2 \; = \;
 [r_{l-1} \, + \, \triangle r_l/3, r_{l} \, - \, \triangle r_l/3)
  \, \times \, \mathbb{S}^{d-1},                                                                                                                                                 
\notag
\\
& 
T_l^3 \; = \;
 [r_{l} \, - \, \triangle r_l/3, r_{l})
  \, \times \, \mathbb{S}^{d-1}
\quad ( l = 1, 2, \ldots, L_\lambda ).
\end{align}
\par
A next step establishing a Hybrid type inequality is to
  construct a certain support mappings
  $\widetilde{w}_{\lambda,l}$, $\widehat{w}_{\lambda,l}$
   and ${w}_{\lambda,l}$
    by making the best of the function $f_{\lambda}$:
They are given by the solutions of
\allowdisplaybreaks\begin{align}
&
\left\{
\begin{array}{rcl}
&
- \triangle \widetilde{w}_{\lambda,1} \; = \; 0
&
\quad \mathrm{in} \quad [0,{r}_1) \times {\mathbb{S}^{d-1}}
\\
&
\left. \widetilde{w}_{\lambda,1}
 \right\vert_{\{r_1\} \times {\mathbb{S}^{d-1}}}
  \; = \; \left. f_{\lambda}
   \right\vert_{\{{r}_{1}\} \times {\mathbb{S}^{d-1}}},
&
\end{array}
\right.
\label{EQ:Support-1}
\\
&
\left\{
\begin{array}{rcl}
&
- \triangle \widetilde{w}_{\lambda,l} \; = \; 0 
& 
\quad \mathrm{in} \quad T_l
\\
&
\left.
 \widetilde{w}_{\lambda,l} 
  \right\vert_{\{ r_{l-1}\} \times {\mathbb{S}^{d-1}}}
   \; = \;
&
\left. f_{\lambda} \right\vert_{\{ {r}_{l-1}\} \times {\mathbb{S}^{d-1}}}
\\
&
\left.
 \widetilde{w}_{\lambda,l}
  \right\vert_{\{ r_{l}\} \times {\mathbb{S}^{d-1}}}
   \; = \;
&
\left. f_{\lambda} \right\vert_{\{ {r}_{l}\} \times {\mathbb{S}^{d-1}}}
\end{array}
\right.
( l = 2, 3, \ldots, L_\lambda ),
\label{EQ:Support-2}
\\
&
\left\{
\begin{array}{rcl}
&
- \triangle \widehat{w}_{\lambda,1} \; = \; 0 
& 
\quad \mathrm{in} \quad [0,\widehat{r}_1) \times {\mathbb{S}^{d-1}},
\\
&
\left. \widehat{w}_{\lambda,1} 
 \right\vert_{\{\hat{r}_1\} \times {\mathbb{S}^{d-1}}}
  \; = \; \left. \widetilde{w}_{\lambda,2}
   \right\vert_{\{\hat{r}_{1}\} \times {\mathbb{S}^{d-1}}}
&
\end{array}
\right.
\label{EQ:Support-3}
\\
&
\left\{
\begin{array}{rcl}
&
- \triangle \widehat{w}_{\lambda,l} \; = \; 0
&
\quad \mathrm{in} \quad \widehat{T}_l
\\
&
\left.
 \widehat{w}_{\lambda,l}
  \right\vert_{\{ \hat{r}_{l}\} \times {\mathbb{S}^{d-1}}}
   \; = \;
&
\left. \widetilde{\widetilde{w}}_{\lambda,l+1}
 \right\vert_{\{ \hat{r}_{l}\} \times {\mathbb{S}^{d-1}}}
\\
&
\left.
 \widehat{w}_{\lambda,l}
  \right\vert_{\{ \hat{r}_{l-1}\} \times {\mathbb{S}^{d-1}}}
   \; = \;
&
\left. \widetilde{w}_{\lambda,l} \right\vert_{\{ \hat{r}_{l-1}\} \times
{\mathbb{S}^{d-1}}}
\end{array}
\right.
( l = 2, 3, \ldots, L_\lambda -1 ),
\label{EQ:Support-4}
\\
&
\left\{
\begin{array}{rcl}
&
- \triangle \widehat{w}_{\lambda,L_\lambda} \; = \; 0
&
\quad \mathrm{in} \quad \widehat{T}_{L_\lambda}
\\
&
\left.
 \widehat{w}_{\lambda,L_\lambda}
  \right\vert_{\{ \hat{r}_{L_\lambda}\} \times {\mathbb{S}^{d-1}}}
   \; = \;
&
\left. {u}_{\lambda}
 \right\vert_{\{ r \} \times {\mathbb{S}^{d-1}}}
\\
&
\left.
 \widehat{w}_{\lambda,L_\lambda}
  \right\vert_{\{ \hat{r}_{L_\lambda-1}\} \times {\mathbb{S}^{d-1}}}
   \; = \;
&
\left. \widetilde{w}_{\lambda,L_\lambda} 
 \right\vert_{\{ \hat{r}_{L_\lambda-1}\} \times {\mathbb{S}^{d-1}}},
\end{array}
\right.
\label{EQ:Support-5}
\\
&
\left\{
\begin{array}{rcl}
&
{w}_{\lambda,1} \; = \; \widehat{w}_{\lambda,1}
\end{array}
\right.
\quad \mathrm{on} \quad T_1,
\label{EQ:Support-6}
\\
&
\left\{
\begin{array}{rcl}
&
- \triangle {w}_{\lambda,l} \; = \; 0 
&
\quad \mathrm{in} \quad T_l
\\
&
\left.
 {w}_{\lambda,l}
  \right\vert_{\{ r_{l} \} \times {\mathbb{S}^{d-1}}}
   \; = \;
&
 \left. \widehat{\widehat{w}}_{\lambda,l}
  \right\vert_{\{ r_{l} \} \times {\mathbb{S}^{d-1}}}
\\
&
\left.
 {w}_{\lambda,l}
  \right\vert_{\{ r_{l-1} \} \times {\mathbb{S}^{d-1}}}
   \; = \;
&
 \left. \widehat{w}_{\lambda,l-1}
  \right\vert_{\{ r_{l-1}\} \times {\mathbb{S}^{d-1}}}
\end{array}
\right.
( l = 2, 3, \ldots, L_\lambda -1 ),
\label{EQ:Support-7}
\\
&
\left\{
\begin{array}{rcl}
&
- \triangle {w}_{\lambda,L_\lambda} \; = \; 0 
&
\quad \mathrm{in} \quad T_{L_\lambda}
\\
&
\left.
 {w}_{\lambda,L_\lambda}
  \right\vert_{\{ r_{L_\lambda} \} \times {\mathbb{S}^{d-1}}}
   \; = \;
&
 \left. (v_\lambda - h_{x_0})
  \right\vert_{\{ r \} \times {\mathbb{S}^{d-1}}}
\\
&
\left.
 {w}_{\lambda,L_\lambda}
  \right\vert_{\{ r_{L_\lambda-1} \} \times {\mathbb{S}^{d-1}}}
   \; = \;
&
 \left. \widehat{w}_{\lambda,L_\lambda-1}
  \right\vert_{\{ r_{L_\lambda-1}\} \times {\mathbb{S}^{d-1}}},
\end{array}
\right.
\label{EQ:Support-9}
\end{align}
where
\allowdisplaybreaks\begin{align}
\skew{4}\widehat{\widehat{w}}_{\lambda,l} (r_l,\omega_{d-1}) 
 \; = \; &
\widehat{w}_{\lambda,l} ({r}_l,\omega_{d-1})
\label{EQ:Support-8}
\\
&    
\, - \, ( r_l \, - \, \hat{r}_{l-1})
 \lint_0^1 \, 
  D_\nu \widehat{w}_{\lambda,l}
   (\hat{r}_{l-1} + \tau ( r_l - \hat{r}_{l-1}),\omega_{d-1})
    \, d\tau
\notag
\\
&
\, - \, ( \hat{r}_{l-1} \, - \, {r}_{l-1})
 \lint_0^1 \,
  D_\nu \widehat{w}_{\lambda,l-1}
   (r_{l-1} + \tau ( \hat{r}_{l-1} - r_{l-1}),\omega_{d-1})
    \, d\tau,
\notag
\\
\skew{4}\widetilde{\widetilde{w}}_{\lambda,l} (\hat{r}_{l},\omega_{d-1}) 
 \; = \; &
  \widetilde{w}_{\lambda,l+1} (\hat{r}_{l},\omega_{d-1})
\label{EQ:Support-9} 
\\
&
\, - \, ( \hat{r}_l \, - \, {r}_{l})
 \lint_0^1 \,
  D_\nu \widetilde{w}_{\lambda,l+1}
   (r_{l} + \tau ( \hat{r}_l - r_{l}),\omega_{d-1})
    \, d\tau
\notag
\\
&
\, - \, ( {r}_{l} \, - \, \hat{r}_{l-1})
 \lint_0^1 \,
  D_\nu \widetilde{w}_{\lambda,l}
   (\hat{r}_{l-1} + \tau ( r_{l} - \hat{r}_{l-1}), \omega_{d-1})
    \, d\tau.
\notag
\end{align}
\par
In the following we set 
\allowdisplaybreaks\begin{align}
\rho_l^1 \; = \;
 \hat{r}_{l-1} + \tau ( r_{l} - \hat{r}_{l-1}), \;
\rho_l^2 \; = \;
 r_{l-1} + \tau ( \hat{r}_{l-1} - r_{l-1})
\notag
\end{align}
for any positive number $\tau$ in $(0,1)$.
\par
We state a few fine properties for the mappings $w_{\lambda,l}$,
 $\widehat{w}_{\lambda,l}$ and $\widetilde{w}_{\lambda,l}$ used below:
  To explain it we shall recall a sequence of hyper spherical-harmonics
   $\{ \phi_n^{(\alpha)}\}$ 
    $( n = 0, 1, \ldots ; \alpha = 1, 2, \ldots, N(n))$,
     where hyper spherical-harmonics 
      $\phi_n^{(1)}$, $\phi_n^{(2)}$, $\ldots$, $\phi_n^{(N(n))}$
       are independent components with degree $n$.
\par
Writing $x$ $=$ $(\rho,\omega_{d-1})$, they signify
\allowdisplaybreaks\begin{align}
&  
\widetilde{w}_{\lambda,l} (x)
\; = \; \widetilde{w}_{\lambda,l} (\rho, \omega_{d-1})
\notag
\\
&
\; = \; \sum_{n=1}^\infty \sum_{\alpha=1}^{N(n)}
 \widetilde{a}_n^{(\alpha)} \Bigl( \frac \rho {r_{l}} \Bigr)^n
  \phi_n^{(\alpha)} (\omega_{d-1})
\, + \, 
 \sum_{n=1}^\infty \sum_{\alpha=1}^{N(n)} 
  \widetilde{b}_n^{(\alpha)}
   \Bigl( \frac {r_{l-1}} \rho \Bigr)^{n+d-2}
    \phi_n^{(\alpha)} (\omega_{d-1})
\notag
\\
&
\, + \, \widetilde{a}_0^{(1)}
\notag
\\
\intertext{with}
&
\widetilde{a}_n^{(\alpha)} \; = \;
 \frac {f_{\lambda}^{n, (\alpha)} (t,{r}_{l}) 
        \, - \, 
         f_{\lambda}^{n, (\alpha)} (t,{r}_{l-1}) \tau_l^{n+d-2}}
          {1 - \tau_l^{2n+d-2}},
\notag
\\
&
\widetilde{b}_n^{(\alpha)} \; = \; -
 \frac {f_{\lambda}^{n, (\alpha)} (t,{r}_l) \tau_l^{n}
        \, - \,
        f_{\lambda}^{n, (\alpha)} (t,{r}_{l-1})}
        {1- \tau_l^{2n+d-2}},
\notag
\\
&
\tau_l \; = \; \frac {{r}_{l-1}}{{r}_{l}}, \;
\notag
\\
&
f_{\lambda}^{n, (\alpha)} (t,r_k)
 \; = \; 
  \fint\limits_{{\mathbb{S}^{d-1}}}
   \langle f_{\lambda} (t,r_k,\omega_{d-1}),
    \phi_n^{(\alpha)} (\omega_{d-1}) \rangle
     \, d\omega_{d-1}
\quad ( k = l-1,l),
\notag
\\
&  
\widehat{w}_{\lambda,l} (x)
\; = \; \widehat{w}_{\lambda,l} (\rho, \omega_{d-1})
\notag
\\
&
\; = \; \sum_{n=1}^\infty \sum_{\alpha=1}^{N(n)}
 \widehat{a}_n^{(\alpha)} \Bigl( \frac \rho {\hat{r}_{l}} \Bigr)^n
  \phi_n^{(\alpha)} (\omega_{d-1})
\, + \, 
 \sum_{n=1}^\infty \sum_{\alpha=1}^{N(n)} \widehat{b}_n^{(\alpha)}
  \Bigl( \frac {\hat{r}_{l-1}} \rho \Bigr)^{n+d-2}
   \phi_n^{(\alpha)} (\omega_{d-1})
\notag
\\
&
\, + \, \widehat{a}_0^{(1)}
\notag
\\
\intertext{with}
&
\widehat{a}_n^{(\alpha)} \; = \;
 \frac {\skew{4}\widetilde{\widetilde{w}}_{\lambda,l+1}^{n, (\alpha)} (\hat{r}_{l})
        \, - \, 
         \widetilde{w}_{\lambda,l}^{n, (\alpha)} (\hat{r}_{l-1}) \hat{\tau}_l^{n+d-2}}
          {1 - \hat\tau_l^{2n+d-2}},
\notag
\\
&
\widehat{b}_n^{(\alpha)} \; = \;
 - \frac {\skew{4}\widetilde{\widetilde{w}}_{\lambda,l+1}^{n, (\alpha)} (\hat{r}_l) \hat\tau_l^{n}
         \, - \,
         \widetilde{w}_{\lambda,l}^{n, (\alpha)} (\hat{r}_{l-1})}
         {1 - \hat\tau_l^{2n+d-2}}, \;
\hat\tau_l \; = \; \frac {\hat{r}_{l-1}}{\hat{r}_{l}},
\notag
\\
&
\skew{4}\widetilde{\widetilde{w}}_{\lambda,l+1}^{n, (\alpha)} (\hat{r}_l)
 \; = \; 
   \fint\limits_{{\mathbb{S}^{d-1}}}
    \langle \skew{4}\widetilde{\widetilde{w}}_{\lambda,l+1} (\hat{r}_l,\omega_{d-1}),
     \phi_n^{(\alpha)} (\omega_{d-1}) \rangle
      \, d\omega_{d-1},
\notag
\\
&
\widetilde{w}_{\lambda,l}^{n, (\alpha)} (\hat{r}_{l})
 \; = \;
   \fint\limits_{{\mathbb{S}^{d-1}}}
    \langle \widetilde{w}_{\lambda,l} (\hat{r}_{l-1},\omega_{d-1}),
     \phi_n^{(\alpha)} (\omega_{d-1}) \rangle
      \, d\omega_{d-1},
\notag
\\
&  
w_{\lambda,l} (x)
\; = \; w_{\lambda,l} (\rho, \omega_{d-1})
\notag
\\
&
\; = \; \sum_{n=1}^\infty \sum_{\alpha=1}^{N(n)}
 a_n^{(\alpha)} \Bigl( \frac \rho {r_{l}} \Bigr)^n
  \phi_n^{(\alpha)} (\omega_{d-1})
\, + \, 
 \sum_{n=1}^\infty \sum_{\alpha=1}^{N(n)} b_n^{(\alpha)}
  \Bigl( \frac {r_{l-1}} \rho \Bigr)^{n+d-2}
   \phi_n^{(\alpha)} (\omega_{d-1})
\notag
\\
&
\, + \, a_0^{(1)}
\notag
\\
\intertext{with}
&
a_n^{(\alpha)} \; = \;
 \frac {\skew{4}\widehat{\widehat{w}}_{\lambda,l}^{n, (\alpha)} ({r}_{l})
        \, - \, 
         \widehat{w}_{\lambda,l-1}^{n, (\alpha)} ({r}_{l-1}) \tau_l^{n+d-2}}
          {1 - \tau_l^{2n+d-2}},
\notag
\\
&
b_n^{(\alpha)} \; = \; -
 \frac {\skew{4}\widehat{\widehat{w}}_{\lambda,l}^{n, (\alpha)} ({r}_l) \tau_l^{n}
        \, - \,
         \widehat{w}_{\lambda,l-1}^{n, (\alpha)} ({r}_{l-1})}
          {1 - \tau_l^{2n+d-2}}, \;
\tau_l \; = \; \frac {r_{l-1}}{r_{l}},
\notag
\\
&
\skew{4}\widehat{\widehat{w}}_{\lambda,l}^{n, (\alpha)} (r_l)
 \; = \; 
   \fint\limits_{{\mathbb{S}^{d-1}}}
    \langle \skew{4}\widehat{\widehat{w}}_{\lambda,l} (r_l,\omega_{d-1}),
     \phi_n^{(\alpha)} (\omega_{d-1}) \rangle
      \, d\omega_{d-1},
\notag
\\
&
\widehat{w}_{\lambda,l-1}^{n, (\alpha)} (r_{l-1})
 \; = \;
   \fint\limits_{{\mathbb{S}^{d-1}}}
    \langle \widehat{w}_{\lambda,l-1} (r_{l-1},\omega_{d-1}),
     \phi_n^{(\alpha)} (\omega_{d-1}) \rangle
      \, d\omega_{d-1}.
\notag
\end{align}
\par
We must remark that the $0$-degree Fourier coefficients
 $a_0^{(1)}$, $\widehat{a}_0^{(1)}$ and $\widetilde{a}_0^{(1)}$ are independent of
  $\rho$ by virtue of 
   $- \triangle w_{\lambda,l}$ $=$ $- \triangle \widehat{w}_{\lambda,l}$
    $=$ $- \triangle \widetilde{w}_{\lambda,l}$
     $=$ $0$ and $D_\nu f_\lambda$ $+$ $r$ $\triangle_\tau$
      $f_\lambda$ $=$ $0$: \eqref{EQ:Support}.
       We mention four technical lemmas:
\begin{Lem}[Algebraic Inequalities]\label{LEM:ALG}
For any positive integer $l$ greater than $1$ and less than $L_\lambda$, 
 we have
\allowdisplaybreaks\begin{align}
&
0 \; \le \; 
   \frac {1 \, - \, \tau_{l}^{n+d-2}}
         { 1 \, - \, \tau_{l}^{2n+d-2}}
\, + \, 
 \frac {n \tau_l^{n} (1 \, - \, \tau_{l}^{d-2})}
       { 1 \, - \, \tau_{l}^{2n+d-2}}
\; \le \; C,
\label{INEQ:Alg-1}
\\
&
\biggl\vert
 \frac { 1 \, - \, \tau_{l+1}^{n+d-2}}
       { 1 \, - \, \tau_{l+1}^{2n+d-2}}
  \, - \, \tau_{l}^{n+d-2}
   \frac {1 \, - \, \tau_{l}^{n+d-2}}
         { 1 \, - \, \tau_{l}^{2n+d-2}}
\biggr\vert
\; \le \; C n \triangle \tau_l,
\label{INEQ:Alg-2}
\\
&
\biggl\vert
 \frac {\tau_{l+1}^{n} ( 1 \, - \, \tau_{l+1}^{d-2})}
       { 1 \, - \, \tau_{l+1}^{2n+d-2}}
  \, - \, 
   \frac {\tau_{l}^{2n+d-2} (1 \, - \, \tau_{l}^{d-2})}
         { 1 \, - \, \tau_{l}^{2n+d-2}}
\biggr\vert \, n
\; \le \; C n \triangle \tau_l.
\label{INEQ:Alg-3}
\end{align}
\end{Lem}
\begin{Lem}\label{LEM:Deriv-Support}
For any number $l$ $( l \, = \, 2, 3, \ldots, L_\lambda)$
 and a respective arbitrary point $y$ $=$ $(\rho,\omega_{d-1})$
  in $T_l$ for $w_{\lambda,l}$, $\hat{T}_l$ for
   $\widehat{w}_{\lambda,l}$ and 
    $T_l$ for $\widetilde{w}_{\lambda,l}$,
     it follows that
\allowdisplaybreaks\begin{align}
&
D_\nu w_{\lambda,l} (y)
\; = \; {W}_{\lambda,l}^1 (y) \, + \, {W}_{\lambda,l}^2 (y),
\notag
\\
\intertext{where}
W_{\lambda,l}^1 (\rho,\omega_{d-1})
 \, & = \; 
  \sum_{n=1}^\infty \frac n \rho \sum_{\alpha = 1 }^{N(n)}
   \frac {1 \, - \, \tau_l^{n+d-2}}{1 \, - \, \tau_l^{2n+d-2}}
\Bigl(
 \Bigl( \frac \rho{r_{l}} \Bigr)^{n}
  \, - \, 
   \Bigl( \frac {r_{l-1}}\rho \Bigr)^{n+d-2}
    \Bigr)
     \, \widehat{w}_{\lambda,l-1}^{n,(\alpha)} (r_{l-1})
      \phi_n^{(\alpha)} (\omega_{d-1}),
\notag
\\
{W}_{\lambda,l}^2 (\rho,\omega_{d-1})
 \, & = \; 
  \sum_{n=1}^\infty \frac {n}\rho \sum_{\alpha = 1 }^{N(n)}
   \frac {\tau_l^n (1 \, - \, \tau_l^{d-2})}{1 \, - \, \tau_l^{2n+d-2}}
    \Bigl( \frac {r_{l-1}}\rho \Bigr)^{n+d-2}
     \, \widehat{w}_{\lambda,l-1}^{n,(\alpha)} (r_{l-1})
      \phi_n^{(\alpha)} (\omega_{d-1})
\notag
\\
&
\, - \,
 \frac {d-2}{\rho}
  \sum_{n=1}^\infty \sum_{\alpha = 1 }^{N(n)}
   \frac {1 \, - \, \tau_l^{n}}{1 \, - \, \tau_l^{2n+d-2}}
    \Bigl( \frac {r_{l-1}}\rho \Bigr)^{n+d-2}
     \, \widehat{w}_{\lambda,l-1}^{n,(\alpha)} (r_{l-1}) 
      \phi_n^{(\alpha)} (\omega_{d-1}),
\notag
\\
&
D_\nu \widehat{w}_{\lambda,l} (y)
 \; = \; \widehat{W}_{\lambda,l}^1 (y) \, + \, \widehat{W}_{\lambda,l}^2 (y),
\notag
\\
\intertext{where}
\widehat{W}_{\lambda,l}^1 (\rho,\omega_{d-1}) 
 \, & = \; 
  \sum_{n=1}^\infty \frac n \rho \sum_{\alpha = 1}^{N(n)}
   \frac {1 \, - \, \hat\tau_l^{n+d-2}}{1 \, - \, \hat\tau_l^{2n+d-2}}
\biggl(
 \Bigl( \frac \rho{\hat{r}_{l}} \Bigr)^{n}
  \, - \, 
   \Bigl( \frac {\hat{r}_{l-1}}\rho \Bigr)^{n+d-2}
    \biggr)
     \, \widetilde{w}_{\lambda,l}^{n,(\alpha)} (\hat{r}_{l-1}) 
      \phi_n^{(\alpha)} (\omega_{d-1}),
\notag
\\
\widehat{W}_{\lambda,l}^2 (\rho,\omega_{d-1})
 \, & = \;
  \sum_{n=1}^\infty \frac n \rho \sum_{\alpha = 1}^{N(n)}
   \frac {\hat\tau_l^n (1 \, - \, \hat\tau_l^{d-2})}{1 \, - \, \hat\tau_l^{2n+d-2}}
\Bigl( \frac {\hat{r}_{l-1}}\rho \Bigr)^{n+d-2}
 \, \widetilde{w}_{\lambda,l}^{n,(\alpha)} (\hat{r}_{l-1})
  \phi_n^{(\alpha)} (\omega_{d-1})
\notag
\\
&
\, - \,
 \frac {d-2}\rho
  \sum_{n=1}^\infty \sum_{\alpha = 1 }^{N(n)}
   \frac {1 \, - \, \hat\tau_l^{n}}
         {1 \, - \, \hat\tau_l^{2n+d-2}}
    \Bigl( \frac {\hat{r}_{l-1}}\rho \Bigr)^{n+d-2}
     \widetilde{w}_{\lambda,l}^{n,(\alpha)} (\hat{r}_{l-1})      
      \, \phi_n^{(\alpha)} (\omega_{d-1}),
\notag
\\
D_\nu \widetilde{w}_{\lambda,l} (y)
&
\notag
\\
\; = \; &
 \sum_{n=1}^\infty \sum_{\alpha = 1}^{N(n)} 
  \frac {f_{\lambda}^{n,(\alpha)} (t,r_l)
         \, - \, 
          f_{\lambda}^{n,(\alpha)} (t,r_{l-1})}
         {1 \, - \, {\tau}_l^{2n+d-2}}
\Bigl(
 \frac n\rho \Bigl( \frac \rho {{r}_{l}} \Bigr)^{n}
  \, + \, {\tau}_l^{n}
   \frac {n+d-2}\rho \Bigl( \frac {{r}_{l-1}}\rho \Bigr)^{n+d-2}
    \Bigr)
     \, \phi_n^{(\alpha)} (\omega_{d-1})
\notag
\\
\, + \, &
 \sum_{n=1}^\infty \frac n\rho \sum_{\alpha = 1 }^{N(n)}
  \frac {1 \, - \, {\tau}_l^{n+d-2}}{1 \, - \, {\tau}_l^{2n+d-2}}
\biggl(
 \Bigl( \frac \rho{r_{l}} \Bigr)^{n} 
  \, - \, 
   \Bigl( \frac {r_{l-1}}\rho \Bigr)^{n+d-2}
    \biggr)
\, f_{\lambda}^{n,(\alpha)} (t,r_{l-1}) \phi_n^{(\alpha)} (\omega_{d-1})
\notag
\\
\, + \, &
 \sum_{n=1}^\infty \frac {n}\rho \sum_{\alpha = 1 }^{N(n)}
  \frac {\tau_l^n (1 \, - \, {\tau}_l^{d-2})}{1 \, - \, {\tau}_l^{2n+d-2}}
\Bigl( \frac {r_{l-1}}\rho \Bigr)^{n+d-2}
 \, f_{\lambda}^{n,(\alpha)} (t,r_{l-1}) \phi_n^{(\alpha)} (\omega_{d-1})
\notag
\\
\, - \, & \frac {d-2}\rho
 \sum_{n=1}^\infty \sum_{\alpha = 1}^{N(n)}
  \frac {1 \, - \, {\tau}_l^{n}}{1 \, - \, {\tau}_l^{2n+d-2}}
   \Bigl( \frac {r_{l-1}} \rho \Bigr)^{n+d-2}
\, f_{\lambda}^{n,(\alpha)} (t,r_{l-1}) \phi_n^{(\alpha)} (\omega_{d-1}).
\notag
\end{align}
\end{Lem}
\begin{Lem}\label{LEM:W-1}
For any positive number $\rho$,
 the symbol $S_\rho (\omega_{d-1})$ denotes by
  a geodesic ball on $\mathbb{S}^{d-1}$ 
   centred at $\omega_{d-1}${\rm;}
    Then for any point $y$ $\in$ $T_l^2$,
\allowdisplaybreaks\begin{align}
&
\vert
 D_\nu w_{\lambda,l} (y)
  \vert^2
\; \le \; C \fint\limits_{B_{{\triangle r_l}/6} (y)}
 \vert 
  D_\nu w_{\lambda,l} (y)
   \vert^2 \, dy
\label{INEQ:W-1}
\\
\intertext{and for $y$ $\in$ $T_l^1$ and $T_l^3$,}
&
\vert
 D_\nu w_{\lambda,l} (y)
  \vert^2
\; \le \; 2
 \fint\limits_{B_{{\triangle r_l}/6} (y)}    
  \vert
   D_\nu \widehat{w}_{\lambda,l-1} (x)
    \vert^2 \, dx
\label{INEQ:W-2}
\\
&
\, + \, C 
 \fint\limits_{[r_{l-1}, r_{l-1} + \triangle r_l /2) \times S_{\triangle \tau_l/2} (\omega_{d-1})}
  \Bigl(
   \vert
    D_\nu w_{\lambda,l} (x)
     \vert^2
      \, + \, 
       \vert
        D_\nu \widehat{w}_{\lambda,l-1} (x)
         \vert^2
          \Bigr) \, dx,
\notag
\\
&
\vert
 D_\nu w_{\lambda,l} (y)
  \vert^2
\; \le \; 2
 \fint_{B_{{\triangle r_l}/6} (y)}    
  \vert
   D_\nu \widehat{w}_{\lambda,l} (x)
    \vert^2 \, dx
\label{INEQ:W-3}
\\
&
\, + \, C 
 \fint\limits_{[r_{l} - \triangle r_l /2, r_{l}) \times S_{\triangle \tau_l/2} (\omega_{d-1})}
  \bigl(
   \vert
    D_\nu w_{\lambda,l} (x)
     \vert^2
      \, + \, 
       \vert
        D_\nu \widehat{w}_{\lambda,l} (x)
         \vert^2
          \bigr) \, dx
\notag
\end{align}
respectively holds.
Furthermore \eqref{INEQ:W-1}, \eqref{INEQ:W-2} and \eqref{INEQ:W-3} imply
\allowdisplaybreaks\begin{align}
&
\sum_{l=2}^{L_\lambda-1}
 \frac {r^2}{\triangle r_l^2}
  \lint_{T_l}
   \vert W_{\lambda,l}^1 \vert^2
    \, dx
\, + \, 
 \sum_{l=2}^{L_\lambda-1}
  \frac {r^2}{\triangle r_l^2}
   \lint_{\hat{T}_l}
    \vert \widehat{W}_{\lambda,l}^1 \vert^2
     \, dx
\notag
\\
&
\; \le \; C r 
 \lint_{\partial B_r}
  \vert D_\tau v_\lambda \vert^2 \, d\mathcal{H}_x^{d-1}
\, + \, o(1)
 \quad ( \lambda \, \nearrow \, \infty ),
\label{INEQ:W-4}
\\
&
\vert W_{\lambda,l}^2 (y) \vert 
 \; \le \; \frac Cr, \quad
  \vert \widehat{W}_{\lambda,l}^2 (y) \vert
   \; \le \; \frac Cr
\quad
( l \, = \, 1,2, \ldots, L_\lambda)
\label{INEQ:W-5}
\\
&
\quad \mathit{for} \; \mathit{any} \; \mathit{respective} 
 \; \mathit{point} \;
  y \, \in \, T_l \; \mathit{and} \; \hat{T}_l,
\notag
\\
&
\sum_{l=2}^{L_\lambda-1} \triangle r_l
\lint_{\{ r_l\} \times \mathbb{S}^{d-1}}
 ( \vert W_{\lambda,l}^2 \vert^2
    \, + \,
     \vert \widehat{W}_{\lambda,l}^2 \vert
  )
\, d\mathcal{H}_x^{d-1}
\notag
\\
&
\; \le \; \frac Cr
 \lint_{\partial B_r}
  | v_\lambda \, - \, h_{x_0} |^2 \, d\mathcal{H}_x^{d-1},
\label{INEQ:W-6}
\\
&
\sum_{l=1}^{L_\lambda-1}
 \frac {1}{\triangle r_l}
\lint_{\{ r_l\} \times \mathbb{S}^{d-1}}
 \vert
  D_\nu 
   ( w_{\lambda,l+1}
    \, - \, 
     w_{\lambda,l})
      \vert^2
       \, d\mathcal{H}_x^{d-1}
\label{INEQ:W-7}
\\
&
\; \le \; C r
 \lint_{\partial B_r}
  | D_\tau v_\lambda |^2 \, d\mathcal{H}_x^{d-1}
\, + \,
 \frac C{r}
  \lint_{\partial B_r}
   | v_\lambda \, - \, h_{x_0} |^2 \, d\mathcal{H}_x^{d-1}
\, + \, o(1)
 \quad ( \lambda \, \nearrow \, \infty ),
\notag
\\
&
\sum_{l=2}^{L_\lambda} \triangle r_l
 \lint_{\{ {\rho_l^1 \} \times \mathbb{S}^{d-1}}}
\vert
 D_\tau D_\nu \widehat{w}_{\lambda,l}
  \vert^2
   \, d\mathcal{H}_x^{d-1}
\notag                                 
\\
&
\; \le \; C \epsilon_0^4 r
 \lint_{\partial B_r}
  \vert D ( v_{\lambda} \, - \, h_{x_0} ) \vert^2
   \, d\mathcal{H}_x^{d-1}
\, + \, \frac {C(\epsilon_0)}r
 \lint_{\partial B_r}
  \vert v_{\lambda} \, - \, h_{x_0} \vert^2
   \, d\mathcal{H}_x^{d-1}
\notag
\\
&
\, + \, o(1)
 \quad ( \lambda \, \nearrow \, \infty ).
\label{INEQ:W-8}
\end{align}
\end{Lem}
\par
From K.Horihata~\cite[Lemma 2.3]{horihata}, we deduce \eqref{INEQ:W-8}.
 The final lemma is as follows: 
\begin{Lem}\label{LEM:Higher-Deriv}
For any positive number $r$, we have
\allowdisplaybreaks\begin{align}
& r^2
\sum_{l=2}^{L_\lambda}
 \lint_{T_l}
  \bigl(
   \vert
    D_\nu^2 \widehat{w}_{\lambda,l}
     \vert^2
\, + \, 
 \vert 
  D_\nu D_\tau \widehat{w}_{\lambda,l}
   \vert^2
\, + \,
 \vert D_\tau^2 \widehat{w}_{\lambda,l}
  \vert^2 \bigr) \, dy
\notag
\\
&
\; \le \; Cr \lint_{\partial B_r}
 \vert D_\tau ( v_{\lambda} \, - \, h_{x_0} )
  \vert^2 \, d\mathcal{H}_y^{d-1}
\, + \, o(1)
 \quad ( \lambda \nearrow \infty)
\label{INEQ:Higher-Deriv-1}
\intertext{and}
& r^2
\lint_{T_1}
 \bigl(
  \vert
   D_\nu^2 \widehat{w}_{\lambda,1}
    \vert^2
\, + \,
 \vert
  D_\nu D_\tau \widehat{w}_{\lambda,1}
   \vert^2
\, + \,
 \vert D_\tau^2 \widehat{w}_{\lambda,1}
  \vert^2 \bigr) \, dy
\notag
\\
&
\; \le \; \frac Cr \lint_{\partial B_r}
 \vert v_{\lambda} \, - \, h_{x_0}                                                            
  \vert^2 \, d\mathcal{H}_y^{d-1}.
\label{INEQ:Higher-Deriv-2}
\end{align}
\end{Lem}
%
%
\vskip 2pt
\par
After the preparation above, 
 we demonstrate the proof of Theorem \ref{THM:HI}.
\vskip 9pt
\noindent{\underbar{Proof of Theorem \ref{THM:HI}.}}
\vskip 9pt
\rm\enspace
Take the difference between \eqref{EQ:GLHF} and 
 $\triangle w_{\lambda,l}$ $=$ $0$ on $T_l$,
  multiplying it by 
   $-2y \cdot D ( (v_\lambda - h_{x_0}) \, - \, w_{\lambda,l})$,
    integrate it on $T_l$ and sum up it for $l$ to verify
\allowdisplaybreaks\begin{align}
&
-2 \, \sum_{l=1}^{L_\lambda} 
 \lint_{T_l}
\left\langle
 \frac {\partial v_\lambda}{\partial t},
  y \!\cdot\! D ( ( v_\lambda - h_{x_0} ) \, - \, w_{\lambda,l} )
 \right\rangle \,  \, dy
\notag
\\
&
\, + \, (d-2) \, 
 \sum_{l=1}^{L_\lambda} 
   \lint_{T_l}
   | D ( ( v_\lambda - h_{x_0} ) \, - \, w_{\lambda,l} )|^2 \,  \, dy
\, + \, \frac {d \, \lambda^{1-\kappa}}{2} \,
 \sum_{l=1}^{L_\lambda} 
   \lint_{T_l}
   ( | v_\lambda |^2 \, - \, 1 )^2 \,  \, dy
\notag
\\
&
\, + \, 2 \sum_{l=1}^{L_\lambda} \lint_{T_l}
\langle
 a\!\cdot\!D ( ( v_\lambda - h_{x_0} ) \, - \, w_{\lambda,l} ),
  D_d ( ( v_\lambda - h_{x_0} ) \, - \, w_{\lambda,l} )
  \rangle \, dy
\notag
\\
&
\, + \, 2 \sum_{i,j=1}^{d-1} \sum_{l=1}^{L_\lambda}
 \lint_{T_l}
  y^j
   \Bigl( \frac {\partial a_j}{\partial y_i}
    \, - \, \frac {\partial a_i}{\partial y_j} \Bigr)
\langle
 D_i ( ( v_\lambda - h_{x_0} ) \, - \, w_{\lambda,l} ),
  D_d ( ( v_\lambda - h_{x_0} ) \, - \, w_{\lambda,l} )
  \rangle \, dy
\notag
\\
&
\; = \; 2 \lambda^{1-\kappa}
 \sum_{l=1}^{L_\lambda} 
   \lint_{T_l}
   ( | v_\lambda |^2 \, - \, 1 )
    \la v_\lambda, y\!\cdot\!D h_{x_0} \ra
     \, dy
\notag\\
&
\, + \, 2 \lambda^{1-\kappa}
 \sum_{l=1}^{L_\lambda}
   \lint_{T_l}
   ( | v_\lambda |^2 \, - \, 1 )
    \la v_\lambda, y \cdot D w_{\lambda,l} \ra
     \, dy
\notag\\
&
\, - 2 \,
 \sum_{l=1}^{L_\lambda}
  \lint_{T_l}
\la ( L \, - \, \triangle ) w_{\lambda,l},
 y\!\cdot\!D ( ( v_\lambda - h_{x_0} ) \, - \, w_{\lambda,l})
  \ra \, dy
\notag
\\
&
\, - \,
 \sum_{l=1}^{L_\lambda} 
  \,
\Bigl( \lint_{ \{ r_{l} \} \times {\mathbb{S}^{d-1}}}
 \, - \,
  \lint_{ \{ r_{l-1} \} \times {\mathbb{S}^{d-1}}}
   \Bigr)
\, 
 \Bigl(
  \rho \, - \, \frac {y_d a\!\cdot\!y}{\sqrt{\rho^2 - y_d^2}}
   \Bigr) \, 
 | D_{\nu} ( ( v_\lambda - h_{x_0} ) \, - \, w_{\lambda,l})|^2
\, d\mathcal{H}_y^{d-1}
\notag\\
&
\, + \,
 \sum_{l=1}^{L_\lambda}
  \,
\Bigl( 
 \lint_{ \{ r_{l} \} \times {\mathbb{S}^{d-1}}}
  \, - \,
   \lint_{ \{ r_{l-1} \} \times {\mathbb{S}^{d-1}}}
    \Bigr)
\,
 \Bigl(
  \rho \, - \, \frac {y_d a\!\cdot\!y}{\sqrt{\rho^2 - y_d^2}}
   \Bigr) \,
| D_{\tau} ( ( v_\lambda - h_{x_0} ) \, - \, w_{\lambda,l})|^2
 \, d\mathcal{H}_y^{d-1}
\notag
\\
&
\, + \,
 \sum_{l=1}^{L_\lambda}
  \,
\Bigl( 
 \lint_{ \{ r_{l} \} \times {\mathbb{S}^{d-1}}}
  \, - \,
   \lint_{ \{ r_{l-1} \} \times {\mathbb{S}^{d-1}}}
    \Bigr)
\notag
\\
& \qquad \qquad \qquad
 \, \frac {\rho y_d}{\sqrt{\rho^2 - y_d^2}}
  \la a\!\cdot\!D_\tau ( ( v_\lambda - h_{x_0} ) \, - \, w_{\lambda,l}),
   D_\nu ( ( v_\lambda - h_{x_0} ) \, - \, w_{\lambda,l})
    \ra \, d\mathcal{H}_y^{d-1}
\notag
\\
&
\, + \,  \frac {\lambda^{1-\kappa}}2
 \sum_{l=1}^{L_\lambda} 
\Bigl(
 \lint_{ \{ r_{l} \} \times {\mathbb{S}^{d-1}}}
  \, - \,
   \lint_{ \{ r_{l-1} \} \times {\mathbb{S}^{d-1}}}
    \Bigr)
\,
 \Bigl(
  \rho \, - \, \frac {y_d a\!\cdot\!y}{\sqrt{\rho^2 - y_d^2}}
   \Bigr) \,
( | v_\lambda |^2 \, - \, 1 )^2 \, d\mathcal{H}_y^{d-1}
\notag
\\
&
\; = \; \; (\mathrm{\bigroman{1}}) \, + \, (\mathrm{\bigroman{2}}) \, + \,
 \, (\mathrm{\bigroman{3}}) \, + \, \cdots
  \, + \,  (\mathrm{\bigroman{7}})
\label{EQ:RP-0}
\end{align}
with
\begin{math}
a\!\cdot\!y \; = \; \sum_{i=1}^{d-1}
 a_i y_i, \;
a\!\cdot\!D \; = \; \sum_{i=1}^{d-1}
 a_i D_i.
\end{math}
\par
From now on  we shall estimate the each term of 
 the right-hand side in \eqref{EQ:RP-0}.
Primarily we estimate ({\bigroman{1}}): 
 By using the usual Sobolev inequality, we infer
\allowdisplaybreaks\begin{align}
(\mathrm{\bigroman{1}}) & \, \le \; 
 \epsilon_0^2 \, \lambda^{1-\kappa}
  \lint_{B_r}
   ( 1 \, - \, | v_\lambda |^2 ) \, dy \;
\notag
\\
&
\, + \, \frac {C r^2 \lambda^{1-\kappa}}{\epsilon_0^2}
 \lint_{B_r}
  ( 1 \, - \, | v_\lambda |^2 ) \, dy \;
\fint\limits_{B_{2r}}
 ( \vert D^{[d/2]+2} h_{x_0} \vert^2 
  \, + \, 
   \vert D h_{x_0} \vert^2
    )
     \, dy.
\label{INEQ:L-1}
\end{align}
\par
Next we estimate ({\bigroman{2}}):
 We choose a collection of balls
  $\{ B_{\triangle r_1/12} (y_i )\}$ $(i \, \in \, I_1)$ with
\allowdisplaybreaks\begin{align}
&  
B_{r_1} \, \subset \, \bigcup_{i \in I_1}
 B_{\triangle r_1/12} (y_i)
  \, \subset \, 
   B_{\hat{r}_1}
\notag
\end{align}
and for any integer $l$ $(l \, = \, 2,3,\ldots,L_\lambda)$,
 do three ones of balls
  $\{ B_{\triangle r_l/12} (y_i )\}$
   $(i \, \in \, I_l^1)$,
    $\{ B_{\triangle r_l/12} (y_i )\}$
     $(i \, \in \, I_l^2)$ and
      $\{ B_{\triangle r_l/12} (y_i )\}$
       $(i \, \in \, I_l^3)$
\allowdisplaybreaks\begin{align}
&  
T_l^1 \, \subset \, \bigcup_{i \in I_l^1}
 B_{\triangle r_l/12} (y_i),
  \quad
T_l^2 \, \subset \, \bigcup_{i \in I_l^2}
 B_{\triangle r_l/12} (y_i),
  \quad
T_l^3 \, \subset \, \bigcup_{i \in I_l^3}
 B_{\triangle r_l/12} (y_i),
  \quad
\notag
\\
&
\max_{j \, \in \, I_l^1 \cup I_l^2 \cup I_l^3}
 \mathrm{Card}
  \{ i \, \in \, I_l^1 \cup I_l^2 \cup I_l^3 \, ; \, 
   B_{\triangle r_l/6} (y_i)
    \cap 
     B_{\triangle r_l/6} (y_j)
      \; \ne \; \emptyset
       \} 
\notag
\end{align}
is finite and independent of $r$ and $l$.
From Lemma \ref{LEM:GLQE}, we obtain
\allowdisplaybreaks\begin{align}
&
\lambda^{1-\kappa}
 \lint_{B_{\triangle r_l/12}(y_i)}
   ( 1 \, - \, | v_\lambda |^2 ) \, dy
\; \le \; C 
 \lint_{B_{\triangle r_l/6}(y_i)}
  \mathbf{e} \, dy
\, + \, \frac C {\triangle r_l} 
 \lint_{B_{\triangle r_l/6}(y_i)}
  | D v_\lambda | \, dy
\notag
\\
&
\, + \, 
 \lint_{B_{\triangle r_l/6}(y_i)}
  \left\vert 
   \frac {\partial v_\lambda}{\partial t}
    \right\vert \, dy.
\notag
\end{align}
Decompose $D_\nu w_{\lambda,l}$
 into $W_{\lambda,l}^1$ and $W_{\lambda,l}^2$
  in Lemma \ref{LEM:Deriv-Support} and
   recall \eqref{INEQ:W-1}, \eqref{INEQ:W-2} and \eqref{INEQ:W-3}
    to yield
\allowdisplaybreaks\begin{align}
(\mathrm{\bigroman{2}}) & \; \le \; 
 \lambda^{1-\kappa} \sum_{i \in I_1}
\lint_{B_{\triangle r_1/12}(y_i)}
 ( 1 \, - \, | v_\lambda |^2 )
  \, \ | y \cdot D w_{\lambda,1} | \, dy
\notag
\\
&
\, + \, 
 \lambda^{1-\kappa}
  \sum_{l=2}^{L_\lambda-1} \sum_{i \in I_l^1 \cup I_l^2 \cup I_l^3}
\lint_{B_{\triangle r_l/12}(y_i)}
 ( 1 \, - \, | v_\lambda |^2 ) 
  \, \ 
   | y \cdot D w_{\lambda,l} | 
    \, dy
\notag
\\
&
\, + \, 
 \sum_{i_{L_\lambda} \in I_{L_\lambda}^1 \cup I_{L_\lambda}^2 \cup I_{L_\lambda}^3} 
\lint_{B_{\triangle r_{L_\lambda}/12}(y_i) \, \cap \, B_r}
 ( 1 \, - \, | v_\lambda |^2 )
  \, \ 
| y \cdot D {w}_{\lambda,{L_\lambda}} | \, dy
\notag
\\
&
\; \le \; r
 \sum_{i \in I_1}
  \lint_{B_{\triangle r_1/6} (y_i)}
   \mathbf{e}_\lambda
    \, dy \,
\, \fint\limits_{B_{\triangle r_l /6}(y_i)}
 | D_\nu w_{\lambda,1} | \, dy
\notag
\\
&
\, + \, r
 \sum_{i \in I_1}
  \lint_{B_{\triangle r_1/6} (y_i)}
   \Bigl\vert \frac {\partial v_\lambda}{\partial t}
    \Bigr\vert
\, \fint\limits_{B_{\triangle r_l /6}(y_i)}
 | D_\nu w_{\lambda,1} | \, dy
\notag
\\
&
\, + \,
 \sum_{i \in I_1}
  \lint_{B_{\triangle r_1/6} (y_i)}
   \vert D v_\lambda \vert
    \, \frac{dy}{\triangle r_1}
\, \fint\limits_{B_{\triangle r_1/6}(y_i)}
 | D_\nu w_{\lambda,1} | \, dy
\notag
\\
&
\, + \, Cr
 \sum_{l=2}^{L_\lambda-1} \sum_{i \in I_l^1}
  \lint_{B_{\triangle r_l /6} (y_i)}
   \mathbf{e}_\lambda
    \, dy \,
\, \fint\limits_{B_{\triangle r_l /6}(y_i)}
 (
  | W_{\lambda,{l}}^1 |
   \, + \,
    | \widehat{W}_{\lambda,{l-1}}^1 |
     )
      \, dy
\notag
\\
&
\, + \, Cr
 \sum_{l=2}^{L_\lambda-1} \sum_{i \in I_l^1}
  \lint_{B_{\triangle r_l /6} (y_i)}
   \Bigl\vert \frac {\partial v_\lambda}{\partial t}
    \Bigr\vert
     \, dy \,
\fint\limits_{B_{\triangle r_l /6}(y_i)}
 (
  | W_{\lambda,{l}}^1 |
   \, + \,
    | \widehat{W}_{\lambda,{l-1}}^1 |
     )
      \, dy
\notag
\\
&
\, + \, C r
 \sum_{l=2}^{L_\lambda-1} \sum_{i \in I_l^1}
  \lint_{B_{\triangle r_l /6} (y_i)}
   \vert D v_\lambda \vert
    \, \frac{dy}{\triangle r_l}
\fint\limits_{B_{\triangle r_l /6}(y_i)}
 (
  | W_{\lambda,{l}}^1 |
   \, + \,
    | \widehat{W}_{\lambda,{l-1}}^1 |
     )
      \, dy
\notag
\\
&
\, + \, Cr
 \sum_{l=2}^{L_\lambda-1} \sum_{i \in I_l^2}
  \lint_{B_{\triangle r_l/6} (y_i)}
   \mathbf{e}_\lambda \, dy \,
\fint\limits_{B_{\triangle r_l /6}(y_i)}
 | W_{\lambda,{l}}^1 | \, dy
\notag
\\
&
\, + \, r
 \sum_{l=2}^{L_\lambda-1} \sum_{i \in I_l^2}
  \lint_{B_{\triangle r_l /6} (y_i)}
   \Bigl\vert \frac {\partial v_\lambda}{\partial t}
    \Bigr\vert
     \, dy \,
\fint\limits_{B_{\triangle r_l /6}(y_i)}
 | W_{\lambda,{l}}^1 | \, dy
\notag
\\
&
\, + \, Cr
 \sum_{l=2}^{L_\lambda-1} \sum_{i \in I_l^2}
  \lint_{B_{\triangle r_l/6} (y_i)}
   \vert D v_\lambda \vert
    \, \frac{dy}{\triangle r_l}
\fint\limits_{B_{\triangle r_l/6}(y_i)}
 | W_{\lambda,{l}}^1 |
  \, dy
\notag
\\
&
\, + \, C r
 \sum_{l=2}^{L_\lambda-1} \sum_{i \in I_l^3}
  \lint_{B_{\triangle r_l/6} (y_i)}
   \mathbf{e}_\lambda \, dy \,
\fint\limits_{B_{\triangle r_l/6}(y_i)}
 (
  | W_{\lambda,l}^1 |
   \, + \,
    | \widehat{W}_{\lambda,l}^1 |
     )
      \, dy
\notag
\\
&
\, + \, r
 \sum_{l=2}^{L_\lambda-1} \sum_{i \in I_l^3}
  \lint_{B_{\triangle r_l/6} (y_i)}
   \Bigl\vert \frac {\partial v_\lambda}{\partial t}
    \Bigr\vert
     \, dy \,
\fint\limits_{B_{\triangle r_l/6}(y_i)}
 (
  | W_{\lambda,l}^1 |
   \, + \,
    | \widehat{W}_{\lambda,l}^1 |
     )
      \, dy
\notag
\\
&
\, + \, C r
 \sum_{l=2}^{L_\lambda-1} \sum_{i \in I_l^3}
  \lint_{B_{\triangle r_l/6} (y_i)}
   \vert D v_\lambda \vert
    \, \frac{dy}{\triangle r_l}
\fint\limits_{B_{\triangle r_l/6}(y_i)}
 (
  | W_{\lambda,l}^1 |
   \, + \,
    | \widehat{W}_{\lambda,l}^1 |
     )
      \, dy
\notag
\\
&
\, + \, r \lambda^{1-\kappa}
 \sum_{l=2}^{L_\lambda-1} \sum_{i \in I_l^1}
  \lint_{B_{\triangle r_l/6} (y_i)}
   ( 1 \, - \, | v_\lambda|^2 )
    \, dy
\fint\limits_{B_{\triangle r_l/6}(y_i)}
 (
  | W_{\lambda,{l}}^2 |
   \, + \,
    | \widehat{W}_{\lambda,{l-1}}^2 |
     )
      \, dy
\notag
\\
&
\, + \, r \lambda^{1-\kappa}
 \sum_{l=2}^{L_\lambda-1} \sum_{i \in I_l^2}
  \lint_{B_{\triangle r_l/6} (y_i)}
   ( 1 \, - \, | v_\lambda|^2 )
    \, dy
\fint\limits_{B_{\triangle r_l/6}(y_i)}
 | W_{\lambda,{l}}^2 |
  \, dy
\notag
\\
&
\, + \, r \lambda^{1-\kappa}
 \sum_{l=2}^{L_\lambda-1} \sum_{i \in I_l^3}
  \lint_{B_{\triangle r_l/6} (y_i)}
   ( 1 \, - \, | v_\lambda|^2 )
    \, dy
\fint\limits_{B_{\triangle r_l/6}(y_i)}
 (
  | W_{\lambda,l}^2 |
   \, + \,
    | \widehat{W}_{\lambda,l}^2 |
     )
      \, dy
\notag
\\
&
\, + \, r \lambda^{1-\kappa}
 \sum_{i \in I_{L_\lambda}^1 \cup I_{L_\lambda}^2 \cup I_{L_\lambda}^3}
  \lint_{B_{\triangle r_{L_\lambda}/12}(y_i) \, \cap \, B_r}
   ( 1 \, - \, | v_\lambda |^2 )
\, | y \cdot D {w}_{\lambda,{L_\lambda}} | \, dy
\notag
\\
&
\; \le \; \epsilon_0^2
 \sum_{i \in I_1}
  \lint_{B_{\triangle r_1 /6}(y_i)}
   \mathbf{e}_\lambda
    \, dy
\, + \, r^2 
 \sum_{i \in I_1}
  \lint_{B_{\triangle r_1 /6}(y_i)}
   \Bigl\vert \frac {\partial v_\lambda}{\partial t} \Bigr\vert^2
    \, dy
\, + \, \epsilon_0^2
 \sum_{i \in I_1}
  \lint_{B_{\triangle r_1 /6}(y_i)}
   \vert D v_\lambda \vert^2 \, dy
\notag
\\
&
\, + \,
 \frac C {\epsilon_0^2}
  \sum_{i \in I_1}
   \lint_{B_{\triangle r_1 /6}(y_i)}
    \mathbf{e}_\lambda \, dy
\fint\limits_{B_{\triangle r_1 /6}(y_i)}
 | y \cdot D w_{\lambda,1} |^2 \, dy
\notag
\\
&
\, + \, \frac 1{\triangle r_1^2}
 \sum_{i \in I_1}
  \lint\limits_{B_{\triangle r_l/6}(y_i)}
   (
    \vert W_{\lambda,1}^1 \vert^2
     \, + \,
      \vert \widehat{W}_{\lambda,1}^1 \vert^2
       )
        \, dy
\notag
\\
&
\, + \, \frac C {\epsilon_0^2 \triangle r_1^2}
 \sum_{i \in I_1}
  \lint_{B_{\triangle r_1 /6}(y_i)}
   | y \cdot D w_{\lambda,1} |^2
    \, dy
\notag
\\
&
\, + \, \frac C{\epsilon_0^2}
 \sum_{l=2}^{L_\lambda-1} \sum_{i \, \in \, I_l^1 \cup I_l^2 \cup I_l^3}
  \lint_{B_{\triangle r_l/6}(y_i)}
   \mathbf{e}_\lambda \, dy \,
\, + \, \frac {Cr^2}{\epsilon_0^2}
 \sum_{l=2}^{L_\lambda-1} \sum_{i \, \in \, I_l^1 \cup I_l^2 \cup I_l^3}
  \lint_{B_{\triangle r_l/6}(y_i)}
   \Bigl\vert
    \frac {\partial v_\lambda}{\partial t}
     \Bigr\vert^2 \, dy
\notag
\\
&
\, + \, \frac 1{\epsilon_0^2\triangle r_l^2}
 \sum_{i \in I_1}
  \lint\limits_{B_{\triangle r_l/6}(y_i)}
   (
    \vert W_{\lambda,1}^1 \vert^2
     \, + \,
      \vert \widehat{W}_{\lambda,1}^1 \vert^2
       )
        \, dy
\notag
\\
&
\, + \, C \epsilon_0^2 r^2 
 \sum_{l=2}^{L_\lambda-1} \sum_{i \, \in \, I_l^1}
  \lint_{B_{\triangle r_l/6}(y_i)}
   \mathbf{e}_\lambda
    \, dy \,
\fint\limits_{B_{\triangle r_l/6}(y_i)}
 ( | W_{\lambda,l}^1 |^2
  \, + \, 
   | W_{\lambda,l-1}^1 |^2
    )  \, dy
\notag
\\
&
\, + \, C \epsilon_0^2 r^2
 \sum_{l=2}^{L_\lambda-1} \sum_{i \, \in \, I_l^2}
  \lint_{B_{\triangle r_l/6}(y_i)}
   \mathbf{e}_\lambda
    \, dy \,
\fint\limits_{B_{\triangle r_l/6}(y_i)}
 | W_{\lambda,l}^1 |^2
  \, dy
\notag
\\
&
\, + \,  C \epsilon_0^2 r^2
 \sum_{l=2}^{L_\lambda-1} \sum_{i \, \in \, I_l^3}
  \lint_{B_{\triangle r_l/6}(y_i)}
   \mathbf{e}_\lambda
    \, dy \,
\fint\limits_{B_{\triangle r_l/6}(y_i)}
  ( | W_{\lambda,l}^1 |^2
  \, + \,
   | \widehat{W}_{\lambda,l}^1 |^2
    )  \, dy
\notag
\\
&
\, + \, C \epsilon_0^2 
 \sum_{l=2}^{L_\lambda-1} \sum_{i \, \in \, I_l^1 }
  \lint_{B_{\triangle r_l/6}(y_i)}
   ( | W_{\lambda,l}^1 |^2
    \, + \,
     | \widehat{W}_{\lambda,l-1}^1 |^2
      )  \, dy
\notag
\\
&
\, + \, C \epsilon_0^2 
 \sum_{l=2}^{L_\lambda-1} \sum_{i \, \in \, I_l^2}
  \lint_{B_{\triangle r_l/6}(y_i)}
   | W_{\lambda,l}^1 |^2
    \, dy
\notag
\\
&
\, + \, C \epsilon_0^2
 \sum_{l=2}^{L_\lambda-1} \sum_{i \, \in \, I_l^3}
  \lint_{B_{\triangle r_l/6}(y_i)}
   ( | W_{\lambda,l}^1 |^2
    \, + \,
     | \widehat{W}_{\lambda,l}^1 |^2
      )  \, dy
\notag
\\
&
\, + \, C \epsilon_0^2 
 \sum_{l=2}^{L_\lambda-1} 
  \bigl( \frac r {\triangle r_l} \bigr)^2
   \sum_{i \, \in \, I_l^1 }
    \lint_{B_{\triangle r_l/6}(y_i)}
     ( | W_{\lambda,l}^1 |^2
      \, + \,
       | \widehat{W}_{\lambda,l-1}^1 |^2
        )  \, dy
\notag
\\
&
\, + \, C \epsilon_0^2 
 \sum_{l=2}^{L_\lambda-1}
  \bigl( \frac r {\triangle r_l} \bigr)^2
   \sum_{i \, \in \, I_l^2 }
    \lint_{B_{\triangle r_l/6}(y_i)}
     | W_{\lambda,l}^1 |^2
      \, dy
\notag
\\
&
\, + \, C \epsilon_0^2 
 \sum_{l=2}^{L_\lambda-1}
  \bigl( \frac r {\triangle r_l} \bigr)^2
   \sum_{i \, \in \, I_l^3 }
    \lint_{B_{\triangle r_l/6}(y_i)}
     ( | W_{\lambda,l}^1 |^2
      \, + \,
       | \widehat{W}_{\lambda,l}^1 |^2
        )  \, dy
\notag
\\
&
\, + \, C \lambda^{1-\kappa}
 \sum_{l=2}^{L_\lambda-1} \sum_{i \, \in \, I_l^1 \cup I_l^2 \cup I_l^3}  
  \lint_{B_{\triangle r_l/12}(y_i)}
   ( 1 \, - \, | v_\lambda |^2 )
    \, dy
\notag
\\
&
\, + \, \lambda^{1-\kappa}
 \sum_{i \in I_{L_\lambda}^1 \cup I_{L_\lambda}^2 \cup I_{L_\lambda}^3}
  \lint_{B_{\triangle r_{L_\lambda}/12}(y_i) \, \cap \, B_r}
   ( 1 \, - \, | v_\lambda |^2 )
\, | y \cdot D {w}_{\lambda,{L_\lambda}} | \, dy.
\label{INEQ:GL-E}
\end{align}
On account of Theorem \ref{THM:Mon}, 
\allowdisplaybreaks\begin{align}
&
\frac 1 {(\triangle r_l /6)^{d-2}}
 \lint_{B_{\triangle r_l /6} (y_i)} 
  \mathbf{e}_\lambda (t,y) \, dy
\notag
\\
&
\; \le \;
 \frac C {(\triangle r_l /3)^d}
  \lint_{t- (\triangle r_l /3)^2}^{t}
   \, ds
\lint_{B_{\triangle r_l/3} (y_i)}
 \mathbf{e}_\lambda \, dy
\, + \, C \triangle r_l^2 \sup_{B_{\triangle r_l/3}}
 | D h_{x_0} |^2
\notag
\\
&
\; \le \;
 \lint_{t- (\triangle r_l /3)^2}^t \, ds
\lint_{\Omega}
 \mathbf{e}_\lambda G_{(t+(\triangle r_l /3)^2,x)}
  \, dy
\, + \, C \triangle r_l^2
\notag
\\
&
\; \le \; C
 \lint_{\Omega}
  | D u_0 |^2 \, dy
\, + \, C \triangle r_l^d
\notag
\end{align}
follows for any time $t$ $\in$ $(t_0-R^2, t_0+R^2)$ 
 with a positive number $R$ less than $\sqrt{t_0}/2$
  and any ball $B_{\triangle r_l /6} (y_i)$,
according to \eqref{INEQ:W-3}, \eqref{INEQ:W-4} and \eqref{INEQ:W-5} 
 of Lemma \ref{LEM:W-1},
  which proceeds to our evaluation as follows:
\allowdisplaybreaks\begin{align}
(\mathrm{\bigroman{2}}) & 
\; \le \; C \epsilon_0^2
 \lint_{B_r}
  \mathbf{e}_\lambda
   \, dy
\, + \, \frac {Cr^2} {\epsilon_0^2}
 \lint_{B_r}
  \left\vert \frac {\partial v_\lambda}{\partial t} \right\vert^2 \, dy
\, + \, 
 \lambda^{1-\kappa}
  \lint_{B_r \setminus B_{(1-\epsilon_0^4)r}}
   ( 1 \, - \, | v_\lambda |^2 ) \, dy
\notag
\\
&
\, + \, \frac {C} {\epsilon_0^2}
 \lint_{B_r \setminus B_{(1-\epsilon_0^4)r}}
  \mathbf{e}_\lambda
   \, dy
\, + \, C \epsilon_0^2 r
 \lint_{\partial B_r}
  \vert D ( v_\lambda \, - \, h_{x_0}) \vert^2 \,d\mathcal{H}_y^{d-1}
\notag
\\
&
\, + \, \frac {C(\epsilon_0)} {r}
 \lint_{\partial B_r}
  | v_\lambda - h_{x_0} |^2 \, d\mathcal{H}_y^{d-1}
\, + \, o(1)
 \quad ( \lambda \, \nearrow \, \infty).
\label{INEQ:L-2}
\end{align}
By using Lemma \ref{LEM:Higher-Deriv}
 we similarly asses (\bigroman{3}):
\allowdisplaybreaks\begin{align}
(\mathrm{\bigroman{3}}) &
 \; \le \;
  \frac{\epsilon_0^2}2
\lint_{T_1} 
 \vert D_\nu (( v_\lambda \, - \, h_{x_0}) \, - \, w_{\lambda,1})
  \vert^2 \, dy
\notag
\, + \, \frac {r^2}{2\epsilon_0^2}
 \lint_{T_1}
  \vert ( L - \triangle ) w_{\lambda,1}
   \vert^2 \, dy
\notag
\\
&
\, + \, \frac 1{2\epsilon_0^2}
 \sum_{l=2}^{L_\lambda}
  \lint_{T_l}
   \vert D_\nu (v_\lambda - h_{x_0}) \vert^2 \, dy
\, + \, \frac {(r \epsilon_0)^2}2
 \sum_{l=2}^{L_\lambda}
  \lint_{T_l}
   \vert ( L - \triangle ) w_{\lambda,l}
    \vert^2 \, dy
\notag
\\
&
\; \le \; C \epsilon_0^2
 \lint_{B_r}
  \vert D_\nu ( v_\lambda - h_{x_0} ) \vert^2 \, dy
\, + \, C \epsilon_0^2 r
 \lint_{\partial B_r}
  \vert D_\tau ( v_\lambda - h_{x_0} ) \vert^2 \, d\mathcal{H}_y^{d-1}
\notag
\\
&
\, + \, \frac C{\epsilon_0^2}
 \lint_{B_r \setminus B_{(1-\epsilon_0^4) r}}
  \vert D_\nu ( v_\lambda \, - \, h_{x_0} ) \vert^2
   \, dy
\notag
\\
&
\, + \, \frac {C(\epsilon_0)}{r}
 \lint_{\partial B_r}
  \vert v_\lambda \, - \, h_{x_0} \vert^2
   \, d\mathcal{H}_y^{d-1}
\, + \, o(1)
 \quad ( \lambda \, \nearrow \, \infty).
\label{INEQ:L-3}
\end{align}
\allowdisplaybreaks\begin{align}
(\mathrm{\bigroman{4}}) & 
 \; \le \; 2 \sum_{l=1}^{L_\lambda-1}
  \lint_{ \{ r_l \} \times {\mathbb{S}^{d-1}}}
   \bigl( r_l \, - \, \frac {y_d \, a \cdot y}{\sqrt{r_l^2 - y_d^2}} \bigr)
\langle
 D_{\nu} ( v_\lambda - h_{x_0} ),
  D_{\nu} ( w_{\lambda,l} \, - \, w_{\lambda,l+1} )
   \rangle \, d\mathcal{H}_y^{d-1}
\notag
\\
&
\, - \, 
 \sum_{l=1}^{L_\lambda-1}
  \lint_{ \{ r_l \} \times {\mathbb{S}^{d-1}}}
   \bigl( r_l \, - \, \frac {y_d \, a \cdot y}{\sqrt{r_l^2 - y_d^2}} \bigr)
\langle
 D_{\nu} ( w_{\lambda,l} \, + \, w_{\lambda,l+1}),
  D_{\nu} ( w_{\lambda,l} \, - \, w_{\lambda,l+1})
   \rangle \, d\mathcal{H}_y^{d-1}
\notag
\\
&
\; \le \; \epsilon_0^6 r  
 \lint_{ \{ r_{1} \} \times {\mathbb{S}^{d-1}}}
  | D_{\nu} ( v_\lambda - h_{x_0} ) |^2 \, d\mathcal{H}_y^{d-1}
\notag
\\
&
\, + \, \frac {Cr} {\epsilon_0^6}
 \lint_{ \{ r_{1} \} \times {\mathbb{S}^{d-1}}}
  ( \vert D_{\nu} w_{\lambda,1} \vert^2 
   \, + \,  \vert D_{\nu}  w_{\lambda,2} \vert^2 )
    \, d\mathcal{H}_y^{d-1}
\notag
\\
&
\, + \, \frac C {\epsilon_0^2}
 \sum_{l=2}^{L_\lambda-1} \, \triangle r_l
  \lint_{ \{ r_{l} \} \times {\mathbb{S}^{d-1}}}
   | D_{\nu} ( v_\lambda - h_{x_0} ) |^2 
    \, d\mathcal{H}_y^{d-1}
\notag
\\
&
\, + \, \frac C{\epsilon_0^2}
 \sum_{l=2}^{L_\lambda-1} \, \triangle r_l
  \lint_{ \{ r_{l} \} \times {\mathbb{S}^{d-1}}}
   ( | D_{\nu} w_{\lambda,l} |^2
    \, + \, | D_{\nu} w_{\lambda,l+1} |^2 )
     \, d\mathcal{H}_y^{d-1}
\notag
\\
&
\, + \, \epsilon_0^2
 \sum_{l=2}^{L_\lambda-1} \, \frac {r^2}{\triangle r_l}
  \lint_{ \{ r_{l} \} \times {\mathbb{S}^{d-1}}}
   | D_{\nu} ( w_{\lambda,l} \, - \, w_{\lambda,l+1} )|^2
    \, d\mathcal{H}_y^{d-1}.
\notag
\end{align}
\par
Using Lemma \ref{LEM:W-1}, 
 from a definition on $r_l$, we arrive at
\allowdisplaybreaks\begin{align}
(\mathrm{\bigroman{4}}) & 
 \; \le \; C \epsilon_0^2
  \lint_{B_r} | D ( v_\lambda - h_{x_0} ) |^2 \, dy
\, + \, 
 \frac C {\epsilon_0^2}
  \lint_{B_r \setminus B_{(1-\epsilon_0^4)r}} 
   | D ( v_\lambda - h_{x_0} ) |^2 \, dy
\notag
\\
&
\, + \, C \epsilon_0^2 r
 \lint_{\partial B_r} | D ( v_\lambda - h_{x_0} ) |^2
  \, d\mathcal{H}_y^{d-1}    
\, + \, \frac {C(\epsilon_0)} {r} \lint_{\partial B_r}
 | v_\lambda - h_{x_0} |^2 \, d\mathcal{H}_y^{d-1}
\notag
\\
&
\, + \, o(1)
 \quad ( \lambda \, \nearrow \, \infty).
\label{INEQ:L-4}
\end{align}
\par
From $\skew{4}\widehat{\widehat{w}}_{\lambda,L_\lambda}$
 $=$ $( v_\lambda - h_{x_0} )$ on $\{ r \}$ $\times$ $\mathbb{S}^{d-1}$,
  we obtain
\allowdisplaybreaks\begin{align}
(\mathrm{\bigroman{4}}) & 
 \; = \; 
\sum_{l=1}^{L_\lambda} r_l
 \lint_{\{ r_{l} \} \times {\mathbb{S}^{d-1}}}
  | D_\tau ( ( v_\lambda - h_{x_0} ) \, - \, w_{\lambda, l}) |^2
   \, d\mathcal{H}_y^{d-1}
\notag
\\
&
\, - \, 
 \sum_{l=1}^{L_\lambda} r_{l-1}
  \lint_{\{ r_{l-1} \} \times {\mathbb{S}^{d-1}}}
   | D_\tau ( ( v_\lambda - h_{x_0} ) \, - \, w_{\lambda,l}) |^2
    \, d\mathcal{H}_y^{d-1}
\notag
\\
&
\; = \; 
\sum_{l=1}^{L_\lambda-1} r_l
 \lint_{\{ r_{l} \} \times {\mathbb{S}^{d-1}}}
  | D_\tau ( ( v_\lambda - h_{x_0} ) \, - \, \skew{4}\widehat{\widehat{w}}_{\lambda,l}) |^2
   \, d\mathcal{H}_y^{d-1}
\notag
\\
&
\, - \,
 \sum_{l=1}^{L_\lambda-1} r_{l}
  \lint_{\{ r_{l} \} \times {\mathbb{S}^{d-1}}}
   | D_\tau ( ( v_\lambda - h_{x_0} ) \, - \, \widehat{w}_{\lambda,l}) |^2
    \, d\mathcal{H}_y^{d-1}
\notag
\\
&
\; = \; -2
 \sum_{l=2}^{L_\lambda-1} r_l
  \lint_{\{ r_{l} \} \times {\mathbb{S}^{d-1}}}
   \langle D_\tau ( v_\lambda - h_{x_0} ),
            D_\tau 
             ( \skew{4}\widehat{\widehat{w}}_{\lambda,l}
              \, - \, 
               \widehat{w}_{\lambda,l})
    \rangle
     \, d\mathcal{H}_y^{d-1}
\notag
\\
&
\, + \, 
 \sum_{l=2}^{L_\lambda-1} r_l
  \lint_{\{ r_{l} \} \times {\mathbb{S}^{d-1}}}
   \langle 
    D_\tau ( \skew{4}\widehat{\widehat{w}}_{\lambda,l}
               \, + \,
                \widehat{w}_{\lambda,l}),
                D_\tau
               ( \skew{4}\widehat{\widehat{w}}_{\lambda,l}
              \, - \,
             \widehat{w}_{\lambda,l})
    \rangle
    \, d\mathcal{H}_y^{d-1}.
\notag
\end{align}
\par
By recalling definition on 
 $\skew{4}\widehat{\widehat{w}}_{\lambda,l}$,
  namely using \eqref{EQ:Support-8} and \eqref{INEQ:W-8} in Lemma \ref{LEM:W-1}
   to forego to estimate
\allowdisplaybreaks\begin{align}
(\mathrm{\bigroman{5}}) &
 \; \le \; \frac C{\epsilon_0^2}
\sum_{l=2}^{L_\lambda} \triangle r_{l}
 \lint_{\{ r_{l} \} \times {\mathbb{S}^{d-1}}}
  ( | D_\tau ( v_\lambda - h_{x_0} ) |^2 \, + \, 
   | D_\tau \widehat{w}_{\lambda,l} |^2 ) 
    \, d\mathcal{H}_y^{d-1}
\notag
\\
&
\, + \, C \epsilon_0^2 r^2
 \sum_{l=1}^{L_\lambda-1} \triangle r_{l}
  \lint_{\{ r_{l} \} \times {\mathbb{S}^{d-1}}}
\Bigl( 
 \vert
  D_\tau D_\nu \widehat{w}_{\lambda,l} (\rho_l^1,\omega_{d-1})
   \vert^2
\, + \, 
 \Bigl\vert
  D_\tau D_\nu \widehat{w}_{\lambda,l-1} (\rho_l^2,\omega_{d-1})
   \, \Bigr\vert^2
    \Bigr)
     \, d\mathcal{H}_y^{d-1}
\notag
\\
&
\; \le \; \frac C{\epsilon_0^2}
 \lint_{B_r \setminus B_{(1-\epsilon_0^4)r}} 
  | D ( v_\lambda - h_{x_0} ) |^2 \, dy
\, + \, C \epsilon_0^2 r 
 \lint_{\partial B_r}
  | D ( v_\lambda - h_{x_0} ) |^2
   \, d\mathcal{H}_y^{d-1}
\notag
\\
&
\, + \, \frac {C(\epsilon_0)} r \lint_{\partial B_r}
 \vert v_\lambda - h_{x_0} \vert^2 \, d\mathcal{H}_y^{d-1}
  \, + \, o(1)
   \quad ( \lambda \, \nearrow \, \infty).
\label{INEQ:L-5}
\end{align}
\par
The estimate on (\bigroman{6}) is as same as above.
 Finally we find that $(\mathrm{\bigroman{7}})$ becomes
\allowdisplaybreaks\begin{align}
&
(\mathrm{\bigroman{7}}) \; = \;
 \frac{\lambda^{1-\kappa}}2
  \lint_{\partial B_r}
   \bigl( r - \frac {y_d \, y\!\cdot\!a}{\sqrt{r^2 \, - \, y_d^2}} \bigr)
    ( | v_\lambda |^2 \, - \, 1 )^2 \, d\mathcal{H}_y^{d-1}.
\label{INEQ:L-7}
\end{align}
\par
On the other hand, recalling the smoothness of $a$ given in (B) in page \pageref{SEC:Intro}, 
 we also have the following estimate for the left-hand side in \eqref{EQ:RP-0},
  which is called $(L)$:
\allowdisplaybreaks\begin{align}
\mathrm{(L)} & \, \ge \; \frac {d-2}4
 \lint_{B_{(1-\epsilon_0^4)r}} \mathbf{e}_\lambda \, dy
\, - \, C \, r^2 
 \lint_{B_r}
   \left|
    \frac {\partial v_\lambda}{\partial t}
     \right|^2 \,  \, dy
\label{INEQ:L}
\\
&
\, - \, \frac {C (\epsilon_0)}{r}
 \lint_{\partial B_r}
  | v_\lambda - h_{x_0} |^2 
   \, d\mathcal{H}_y^{d-1}
\, - \, C 
 \lint_{B_r} \vert D h_{x_0} \vert^2 \, dx
\notag
\end{align}
\par
A substitution of \eqref{INEQ:L-1}, \eqref{INEQ:L-2},
 \eqref{INEQ:L-3}, \eqref{INEQ:L-4}, \eqref{INEQ:L-5}, 
  \eqref{INEQ:L-7}, \eqref{INEQ:L} for \eqref{EQ:RP-0}, 
an integration of it with respect to $t$ 
 $\in$ $(-r^2,r^2)$ verifies
\allowdisplaybreaks\begin{align}
\lint_{P_{r/2}} & \mathbf{e}_\lambda \, dz
\; \le \; C \epsilon_0^2
 \lint_{P_r} 
  \mathbf{e}_\lambda \, dz
\, + \, \frac {C r^2}{\epsilon_0^2}
 \lint_{P_{r}} 
  \biggl| \frac {\partial v_\lambda}{\partial t}  \biggr|^2
   \, dz
\, + \,
 \frac {C}{\epsilon_0^2}
  \lint_{-r^2}^{r^2} \, dt \,
   \lint_{B_r \setminus B_{(1-\epsilon_0^4)r}} 
    | D ( v_\lambda - h_{x_0} ) |^2 \, dy
\notag
\\
&
\, + \, C \epsilon_0^2 \lint_{P_r}
 \lambda^{1-\kappa} \,
  ( 1 \, - \, |v_\lambda|^2) \, dz
\notag
\\
&
\, + \, \frac C{\epsilon_0^{10}}
 \lint_{P_r}
  \lambda^{1-\kappa}
   ( 1 \, - \, |v_\lambda|^2) \, dz
\Bigl(
 \fint\limits_{\partial B_{r}}
  \vert v_\lambda \, - \, h_{x_0} \vert^2
   \, d\mathcal{H}_y^{d-1}  
    \, + \, 
     \fint\limits_{\partial B_{r}}
      \vert h_{x_0} \, - \, a \vert^2
       \, d\mathcal{H}_y^{d-1}
        \Bigr)
\notag
\\
&
\, + \, 
 \lint_{-r^2}^{r^2}
  \lambda^{1-\kappa} \, dt
   \lint_{B_r \setminus B_{(1-\epsilon_0^4)r}}
    ( 1 \, - \, |v_\lambda|^2) \, dy
\, + \, \frac{C (\epsilon_0)} {r}
 \lint_{-r^2}^{r^2} \, dt \, 
  \lint_{\partial B_{r}}
   | v_\lambda \, - \, h_{x_0} |^2 \, d\mathcal{H}_y^{d-1}
\notag
\\
&
\, + \, 
 \lint_{-r^2}^{r^2} \frac {r \lambda^{1-\kappa}}2 \, dt  
  \, \lint_{\partial B_{r}}
   ( |v_\lambda|^2 \, - \, 1 )^2 \, d\mathcal{H}_y^{d-1}
\, + \, 
 C \lint_{P_r} \vert D h_{x_0} \vert^2 \, dz
  \, + \, o(1)
   \quad ( \lambda \, \nearrow \, \infty).
\label{INEQ:RP-2}
\end{align}
\par
Integrate \eqref{INEQ:RP-2} from $R/2$ to $R$ with respect to $r$
 and divide it by $R$ to obtain
\allowdisplaybreaks\begin{align}
&
\lint_{P_{R/4}} \mathbf{e}_\lambda \, dz
\; \le \; C \epsilon_0^2
 \lint_{P_{R}} \mathbf{e}_\lambda \, dz
\,+ \,
 \frac {C}{\epsilon_0^2 R}
  \lint_{-R^2}^{R^2} \, dt \,
   \lint_{R/2}^R \, dr
    \lint_{B_r \setminus B_{(1-\epsilon_0^4)r}} | D ( v_\lambda - h_{x_0} ) |^2
     \, dy
\notag
\\
&
\, + \, \frac 1R
 \lint_{-R^2}^{R^2} \lambda^{1-\kappa} \, dt \,
  \lint_{R/2}^R \, dr
   \lint_{B_r \setminus B_{(1-\epsilon_0^4)r}}
    ( 1 \, - \, |v_\lambda|^2) \, dy
\notag
\\
&
\, + \, \frac {C R^2}{\epsilon_0^2}
 \lint_{P_{R}} \,
  \Bigl| \frac {\partial v_\lambda}{\partial t} \Bigr|^2
   \, dz
\, + \,
 \frac {C(\epsilon_0)}{R^2}
  \lint_{P_R} \vert v_\lambda - h_{x_0} \vert^2 \, dz
\notag
\\
&
\, + \, \frac 1{R^{d+2}}
 \lint_{P_R} \lambda^{1-\kappa}
  ( 1 \, - \, |v_\lambda|^2) \, dz
\fint\limits_{P_R}
 \vert v_\lambda \, - \, h_{x_0} \vert^2
  \, dz
\notag
\\
&
\, + \, \frac 1{R^{d-2}}
 \lint_{P_R}
  \lambda^{1-\kappa} 
   ( 1 \, - \, |v_\lambda|^2) \, dz
\fint\limits_{P_R}
  \Bigl\vert \frac {\partial v_\lambda}{\partial t} 
   \Bigr\vert^2 \, dz
\label{INEQ:RP-3}
\\
&
\, + \, 
 \frac 12 
  \lint_{-R^2}^{R^2} \lambda^{1-\kappa} \, dt \,
   \lint_{B_R \setminus B_{R/2}}
    ( | v_\lambda |^2 \, - \, 1 )^2 \, dy
\notag
\\
&
\, + \, C 
 \lint_{P_R} 
  ( \vert D^{[(d+1)/2]+1} h_{x_0} \vert^2 
   \, + \, 
    \vert D h_{x_0} \vert^2 )
     \, dz
\, + \, o(1)
 \quad ( \lambda \, \nearrow \, \infty).
\notag
\end{align}
\par
By interchanging the integration of the second and the third terms
 on the right-hand side in \eqref{INEQ:RP-3} and employing Lemma \ref{LEM:GLQE},
  we claim
\allowdisplaybreaks\begin{align}
\lint_{P_{R/4}} & \mathbf{e}_\lambda \, dz
\; \le \; C \epsilon_0^2
 \lint_{P_{R}} \mathbf{e}_\lambda \, dz
\, + \, \frac {C R^2}{\epsilon_0^2}
 \lint_{P_{R}} \,
  \Bigl| \frac {\partial v_\lambda}{\partial t} \Bigr|^2
   \, dz
\notag
\\
&
\, + \, 
 \frac {C(\epsilon_0)}{R^2} 
  \lint_{P_{R}} | v_\lambda - h_{x_0} |^2 \, dz
\notag
\\
&
\, + \, C
 \lint_{P_R}
  ( \vert D^{[(d+1)/2]+1} h_{x_0} \vert^2
   \, + \,
    \vert D h_{x_0} \vert^2 )
     \, dz
\, + \, o(1)
 \quad ( \lambda \, \nearrow \, \infty).
\label{INEQ:RP-4}
\end{align}
\par
The rest of our proof is perfectly as same as
 the one by K.Horihata~\cite{horihata}.

%
%
\setcounter{chapternumber}{3}\setcounter{equation}{0}
\renewcommand{\theequation}%
           {\thechapternumber.\arabic{equation}}
\section{\enspace WHHF}
\subsection{\enspace Partial Regularity}
This chapter studies a partial boundary regularity on WHHF
 constructed in K.Horihata~\cite{horihata}.
  We review two convergent theorems:
\begin{Thm}{\rm{(Convergence).}}\label{THM:Convergence}
There exist a sub-sequence $\{ u_{\lambda (\nu)} \}$ $( \nu \, = \, 1,2,\ldots )$ of
 $\{ u_\lambda \}$ $( \lambda > 0 )$ in $V_B (Q(T) ; \mathbb{S}^D)$
  and a mapping $u \, \in \,$ $V_B(Q(T) ; \mathbb{S}^D)$                                                                                                                 
   such that the sequence of mappings $\{ u_{\lambda (\nu)} \}$
    $( \nu \, = \, 1,2,\ldots )$
     respectively converges weakly and weakly-$*$ to a mapping $u$ in
      $H^{1,2} (0,T ; L_B^2 (\Omega ; \mathbb{R}^{D+1}))$
       and $L^\infty (0,T ; H_B^{1,2} (\Omega ; \mathbb{R}^{D+1}))$.
So does it strongly to the mapping $u$ in
 $L^2 ((0,T) ; L_B^2 (\Omega ; \mathbb{S}^D))$ and point-wisely to it
  in almost all $z$ $\in$ $Q(T)$ as $\nu \nearrow \infty$.
\end{Thm}
\par
Theorem \ref{THM:Convergence} enables us state the following existence
theorem:
\begin{Thm}{\rm{(Existence).}}\label{THM:Existence-HF}
The GLHF converges to a WHHF in 
 $L^2 ((0,T) ; L_B^2 (\Omega ; \mathbb{S}^D))$
  as $\lambda \nearrow \infty$ \rm{(}\it{modulo a sub-sequence of $\lambda$}\rm{)}.
\end{Thm}
\begin{Def}{}\label{REM:Measure}
Let $\{ u_{\lambda (\nu )}\}$ $( \nu \, = \, 1,2,\ldots )$
 be the sequence selected above and set $\mathbf{e}_{\lambda (\nu)}$
  the Ginzburg-Landau energy density
   $| \nabla u_{\lambda (\nu)} |^2/2$ $+$
    $\lambda (\nu)^{1-\kappa }$ $(| u_{\lambda (\nu)}|^2 \, - \, 1)^2/4$.
We then denotes $\overline{\mathcal{M}}$ by
\begin{equation*}
\overline{\mathcal{M}} (P_R (z_0) \cap Q(T)) \,= \, 
 \limsup_{\lambda (\nu) \nearrow \infty } \frac 1 {R^d}
  \lint_{P_R (z_0) \cap Q(T)} 
   \mathbf{e}_{\lambda (\nu) } \, dz,
\end{equation*}
where $P_R (z_0)$ is an arbitrary parabolic cylinder.
\end{Def}
\begin{Lem}{\rm{(Measured Hybrid Inequality).}}\label{LEM:Measure-Hybrid}
Assume that a sequence of GLHF $\{ u_{\lambda (\nu)} \}$ 
 $( \nu = 1,2,\ldots )$, respectively converges 
  weakly and weakly-$*$ in $H^{1,2} (0,T ; L_B^2 (\mathbb{B}^d ; \mathbb{R}^{D+1}))$  and 
   $L^\infty (0,T;H_B^{1,2} (\mathbb{B}^d ; \mathbb{R}^{D+1}))$
    to a WHHF $u$ $\, \in \,$ $V_B( Q ; \mathbb{S}^D )$ as $\lambda (\nu) \nearrow \infty$.
Then take the pass to the limit $\lambda (\nu) \nearrow \infty$ in 
 Theorem \ref{THM:HI} to infer the following{\rm{:}}
For any positive number $\epsilon_0$,
 there exists a positive constant $C ( \epsilon_0 )$
  satisfying $C ( \epsilon_0 )$ $\nearrow \infty$ as $\epsilon_0  \searrow 0$
   such that the inequality
\allowdisplaybreaks\begin{align}
\overline{\mathcal{M}} & \bigl(P_R (z_0) \cap Q(T) \bigr) 
 \; \le \; 
  \epsilon_0 \overline{\mathcal{M}} \bigl(P_{2R} (z_0) \cap Q(T) \bigr)
\notag
\\
&
\, + \, C ( \epsilon_0 )
 \fint\limits_{P_{2R} (z_0) \cap Q(T)}  | u (z) \, - \, h_0 (x) |^2 \, dz
\notag 
\\
&  
\, + \, C \fint\limits_{P_{2R} \cap Q(T)}
 ( | \nabla^{[(d+1)/2]+1} h_0 (x) |^2 \, + \, | \nabla h_0 (x) |^2 )
  \, dz
\label{INEQ:Measure-Hybrid}
\end{align}
holds for any parabolic cylinder $P_{2R} (z_0)$.
\end{Lem}
Likewise L.Simon~\cite[Lemma 2, p.31]{simon95},
 we can assert the following reverse Poincar\'e inequality.
\begin{Cor}{\rm{(Reverse Poincar\'e Inequality).}}\label{COR:RPI}
The inequality \eqref{INEQ:Measure-Hybrid} implies 
 that the following
\allowdisplaybreaks\begin{align}
& \overline{\mathcal{M}} (P_{R} (z_0) \cap Q(T))
 \; \le \; C 
  \fint\limits_{P_{2R} (z_0) \cap Q(T)} | u (z) \, - \, h_0 (x) |^2 \, dz
\notag
\\
&
\, + \,C \fint\limits_{P_{2R} \cap Q(T)}
 ( | \nabla^{[(d+1)/2]+1} h_0 (x) |^2 \, + \, | \nabla h_0 (x) |^2 ) \, dz
\label{INEQ:RPI}
\end{align}
holds whenever $P_{2R} (z_0)$ is an arbitrary parabolic cylinder.
\end{Cor}
%
By combining Corollary \ref{COR:RPI} with passing to the limit $\lambda$
 $\nearrow$ $\infty$ with Sobolev imbedding theorem and Poincar\'{e} inequality
  for the space variables,
   we can describe the following lemma.
We refer the proof to Theorem 2.1 in M.Giaquinta and M.Struwe~\cite{giaquinta-struwe}.
\begin{Lem}{}\label{LEM:RH}
There exists a positive number $q_0$ greater than $1$
 such that differentials $\nabla u$ of the WHHF $u$ belongs to 
  $L^{2q_0}$ $((0,T) ; L_B^{2q_0} ( \Omega ; \mathbb{R}^{d{(D+1)}}))$
   with
\vskip 15pt
\allowdisplaybreaks\begin{align}
&
\biggl( 
 \fint_{P_{R} (z_0)} | \nabla u (z) |^{2q_0} \, dz
  \biggr)^{1/2q_0}
\; \le \; C
 \biggl( \fint_{P_{2R} (z_0)}
  | \nabla u (z) |^2 \, dz
   \biggr)^{1/2}
\notag
\\
&
\, + \, C
 \biggl( \fint_{P_{2R} (z_0)}
  ( | \nabla^{[(d+1)/2]+1} h_0 (x) |^2 \, + \, | \nabla h_0 (x) |^2 ) \, dz 
   \biggr)^{1/2q_0}
\label{INEQ:RH}
\end{align}
for any parabolic cylinder $P_{2R}$.
\end{Lem}
If we apply \eqref{INEQ:RPI} in Corollary \ref{COR:RPI} to
 the proof on a partial regularity result by Y.Chen~\cite[Lemma 3.1]{chen-91},
  we can claim
\begin{Thm}{}\label{THM:PR}
For any positive number $\epsilon$, set
\begin{align}
&
\mathbf{sing} ( \epsilon ) \, = \, \bigcap_{R > 0}
 \{ z_0 \, \in \, (0,T)\times\overline{\Omega} \, ; \,
  \, \overline{\mathcal{M}} \! \bigl(P_R (z_0)\bigr) \, \ge \, \epsilon \,
   \},
\label{EQ:Sing}
\\
&
\mathbf{reg} ( \epsilon ) \; = \; (0,T)\times\overline{\Omega} \setminus \mathbf{sing}.
\label{EQ:Reg}
\end{align}
Then there exist some positive number $\epsilon_0$ and an increasing
 function $g(t)$ with $g(0) = 0$ and $g(t)$ $=$ $O (t \log (1/t)^{d+1})$
  $( t \searrow 0 )$
   such that
    if $z_0$ $\in$ $reg ( \epsilon_0 )$, that is
     for some positive number $R_0$ and positive integer $\lambda_0$
      possibly depending on $z_0$,
\begin{equation}
\frac 1{R_0^d} \lint_{P_{g(R_0) (z_0)} \cap (0,T)\times\overline{\Omega}}
 \mathbf{e}_\lambda (z) \, dz                                                                                                                           
  \; < \; \epsilon_0                                                                                                                                                              
\end{equation}
implies
\begin{equation}
\sup_{z \in P_{R_0} (z_0) \cap (0,T)\times\overline{\Omega}}                                                                                                                     
 \mathbf{e}_\lambda (z) 
  \; \le \; \frac C{R_0^2}
\label{INEQ:PR}
\end{equation}
as long as any $\lambda$ is more than or equal to $\lambda_0$.
\end{Thm}
\begin{Def}{}
In the sequel we respectively mean $\, \mathbf{sing}$ and $\mathbf{reg}$
  by $\mathbf{sing} ( \epsilon_0 )$ and $\mathbf{reg} ( \epsilon_0 )$.
\end{Def}
\begin{Lem}{}\label{LEM:Cont-Time}
Pick up any point $z_0$ $\in$ $\mathbf{reg}$
 and fix it.
  On the parabolic cylinder $P_{R_0/2} (z_0)$ which is the half size of
   the cylinder in \eqref{INEQ:PR}, we obtain
\begin{equation}
| u_\lambda (t,x) \, - \, u_\lambda (s,x) |
 \; \le \; \frac C{R_0} |t-s|^{1/2}
\label{INEQ:Cont-Time}
\end{equation}
for any points $t$ and $s$ in $[t_0 - (R_0/2)^2, t_0 + (R_0/2)^2]$
 and $x$ $\in$ $\overline{B}_{R_0/2} (x_0)$ with $z_0$ $=$ $(t_0,x_0)$.
\end{Lem}
%
\begin{Thm}{\rm{(Singular Set).}}\label{THM:Singular}
The set of $\mathbf{sing}$ is a closed set and 
\begin{equation}
\mathcal{H}^{(d)} (\mathbf{sing}) \, = \, 0
\label{EQ:Hausdorff-Estimate}
\end{equation}
holds with respect to the parabolic metric.
\end{Thm}
We finally exhibit the strong convergence
 of $\{ u_{\lambda (\nu)} \}$ $( \nu = 1,2, \ldots)$ to a WHHF $u$ 
  in $H^{1,2}$-topology as $\lambda (\nu) \nearrow \infty$.
   The estimates of $\mathcal{H}^{(d)} (\mathbf{sing})$ $\, = \,$ $0$ in Theorem
    \ref{THM:Singular} plays a crucial role in the proof.
However since the proof is 
 as same as in the one in K.Horihata~\cite[Theorem 3.12]{horihata}, 
  we omit it.
\begin{Thm}{\rm{(Strong Convergent of Gradients of \textrm{WHHF}).}}
\label{THM:Strong-Convergence-Gradient}
For a suitable sub-sequence of $\{ \lambda (\nu) \}$
 still denoted by $\{ \lambda (\nu) \}$
  $( \nu = 1,2, \ldots )$, a sequence of the gradients of the GLHF,
   $\{ \nabla u_{\lambda (\nu)}\}$ $( \nu = 1,2, \ldots)$ converges
    strongly to the ones of the WHHF in 
     $L^2 ((0,T) ; L_B^2 (\Omega ; \mathbb{R}^{d(D+1)}))$.
\end{Thm}
%

%
%
\renewcommand{\theequation}%
           {\thechapternumber.\arabic{equation}}
\section{\enspace Proof of Main Theorems}
\par
By utilizing a few ingredients 
 and properties on the WHHF and the GLHF,
  this chapter proves
   Theorem \ref{THM:Main-1}, Theorem \ref{THM:Main-2} 
    and Theorem \ref{THM:Main-3} in Chapter \ref{SEC:Intro}.
\subsection{Proof of Theorem \ref{THM:Main-1}}
We discuss a partial boundary regularity on the WHHF constructed above:
On account of Theorem \ref{THM:Singular} and Theorem \ref{THM:Strong-Convergence-Gradient}, 
 we obtain
\allowdisplaybreaks\begin{align}
&
\mathbf{sing} \, = \, \bigcap_{R > 0}
 \bigl\{ z_0 \in (0,T)\times\overline{\Omega} \, ; \, 
  \frac 1 {2R^d} \kern-2.5pt \lint_{P_R (z_0) \cap Q(T)} \kern-2.5pt
   | \nabla u |^2 \, dz \, \ge \, \epsilon_0 \,
\bigr\},
\label{EQ:Sing-2}
\\
&
\mathbf{reg} \; = \; (0,T)\times\overline{\Omega} \setminus \mathbf{sing},
\label{EQ:Reg-2}
\end{align}
where a number $\epsilon_0$ is a positive constant appeared in Theorem
\ref{THM:PR}.
From Theorem \ref{THM:Singular}, we see that \textbf{sing} is closed. 
The size of \textbf{sing} is analogously measured 
 as in \cite[Theorem 3.10]{horihata};
  So does the smoothness of our WHHF on $\mathbf{reg}$
\par
We readily see that 
 the WHHF $u$ satisfies 
  a monotonical inequality (i) and
   a reverse Poincar\'{e} inequality (ii)
    from Theorem \ref{THM:Mon} with Theorem \ref{THM:Strong-Convergence-Gradient}
     and Corollary \ref{COR:RPI}.
$\qed$
\vskip 6pt
\subsection{Proof of Theorem \ref{THM:Main-2}}
\vskip 6pt
The proof is done by a direct combination of
 Theorem \ref{THM:Mon} and Theorem \ref{THM:Convergence}.
$\qed$
\vskip 6pt
\subsection{Proof of Theorem \ref{THM:Main-3}}
\vskip 6pt
\par
If we recall Theorem \ref{THM:PR}, we have only to show
 that for any point $z_0$ $=$ $(t_0,x_0)$,
\begin{equation}
\lint_{P_r (z_0) \cap Q(T)} | \nabla u |^2 \, dz
 \; < \; \frac {\epsilon_0^2  r^d}2
\end{equation}
holds whereupon $r$ is a positive number.
If we use a stereo-graphic projection 
 given by the mapping 
  $v$ $=$ $( v_i )$ $(i=1,2,\ldots,D)$ of
\begin{equation}
\left\{
\begin{array}{rl}
u^i \, & = \; \displaystyle\frac{2v^i}{1+|v|^2}
\quad ( i \, = \, 1,2,\ldots,D),
\\[4pt]
u^D \, & = \; \displaystyle\frac{1-|v|^2}{1+|v|^2}
\end{array}
\right.
\end{equation}
with $|v|^2$ $=$ $\displaystyle\sum_{i=1}^D$ $(v^i)^2$,
our equation \eqref{EQ:GLHF} becomes for $v^i$ $( i \, = \, 1,2,\ldots,D)$
\begin{equation}
\dfrac {\partial}{\partial t} 
 \bigl( \frac{v^i}{1+|v|^2} \bigr)
  \, = \,
   \triangle \frac{v^i}{1+|v|^2} 
    \, + \, 4 \Bigl| \frac {\nabla v}{1+|v|^2} \Bigr|^2
     \frac{v^i}{1+|v|^2} 
      \quad \mathrm{in} \; Q(T).
\label{EQ:HHF-SG}
\end{equation}
\par
By multiplying \eqref{EQ:HHF-SG} by $v^i$,
 we have
\begin{equation}
\dfrac {\partial W(|v|^2)}{\partial t}
 \, - \, \triangle W(|v|^2)
  \, + \, \frac {4 |\nabla v|^2}{(1+|v|^2)^2}
   \; = \; 0,
\label{EQ:HHF-SG-1}
\end{equation}
where the function $W$ is given by
\begin{equation}
W(x) \; = \; \lint_0^x \frac {(1-t^2) \, dt}{(1+t^2)^2}
 \quad ( x \, \in \, \mathbb{R}).
\notag
\end{equation}
\par
Due to a maximal principle, we verify
\begin{equation}
W(|v|^2) \; \le \; W(1-\theta_0),
\notag
\end{equation}
with $\theta_0$ $=$ 
 $\dist (v_0(\overline{\Omega}), \partial B_1^D)/2$.
\par
While we readily see that there exists a positive integer
 $k_0$ possibly depending on $d$ and $\theta_0$ 
  such that 
\allowdisplaybreaks\begin{align}
&
\esssup_{P_{d_0/2^{k_0-1}} (z_0)} W(|v|^2)
 \, - \,
  \esssup_{P_{d_0/2^{k_0}} (z_0)} W(|v|^2)
\; < \; \epsilon_0^2
\notag
\end{align}
because of 
\allowdisplaybreaks\begin{align}
&      
\sum_{k=1}^\infty
 (\esssup_{P_{d_0/2^{k-1}} (z_0)} W(|v|^2) 
  \, - \, 
   \esssup_{P_{d_0/2^k} (z_0)} W(|v|^2))
\notag
\\
&
\; = \; 
 \esssup_{P_{d_0} (z_0)} W(|v|^2)
  \; \le \; C.
\notag
\end{align}
Consequently, by employing \eqref{EQ:HHF-SG-1} 
 and a sub-linear estimate on $W$,
  we deduce
\begin{equation}
\lint_{P_{d_0/2^k} (z_0) \cap Q(T)}
 | \nabla u |^2 \, dz
  \; < \; C \epsilon_0^2 (d_0/2^k)^d \; \le \; \epsilon_0 (d_0/2^k)^d
\label{INEQ:Small-Energy-Density}
\end{equation}
as long as a positive number $\epsilon_0$ is less than $1/C$.
\par
Moreover from \eqref{INEQ:Small-Energy-Density}, an application on a $\epsilon$-regularity theory in which 
 we exploit \eqref{THM:Existence-HF} or \eqref{INEQ:RPI} to our WHHF $u$ verifies
\begin{equation}
\lint_{P_r (z_0) \cap Q(T)} \vert \nabla u \vert^2 \, dz
 \; < \; C r^{d+\alpha_0}
\end{equation}
for all positive number $r$ less than $d_0/2^{k_0+1}$,
 where $\alpha_0$ is a positive number less than $1$
  and is independent of a point $z_0$ and 
   a positive number $r$.
Thus we arrive at
\begin{equation}
 \lint_{P_{g(r)} (z_0) \cap Q(T)} \vert\nabla u \vert^2 \, dz
  \; < \; \epsilon_0 r^d,
\notag
\end{equation}
where the function $g$ was defined in Theorem \ref{THM:PR}.
\par
Here we selected a positive number $r$ so as to be less than $(C\epsilon_0)^{1/\alpha_0}$,
 which conclude our claim.
$\qed$

%
%
%

\AD

\end{document}